\documentclass[preprint,11pt]{elsarticle}

\newtheorem{theorem}{Theorem}[section]
\newtheorem{corollary}{Corollary}[section]
\newtheorem{lemma}{Lemma}[section]
\newtheorem{remark}{Remark}[section]

\usepackage{amssymb}
\begin{document}

\begin{frontmatter}

\title{Deviation inequalities for martingales with applications to linear regressions and weak invariance principles}
\author{  \ Xiequan Fan }
 \cortext[cor1]{\noindent Corresponding author. \\
\mbox{\ \ \ \ }\textit{E-mail}: fanxiequan@hotmail.com (X. Fan). }
\address {Regularity Team, Inria and MAS Laboratory, Ecole Centrale Paris - Grande Voie des Vignes,\\ 92295 Ch\^{a}tenay-Malabry, France.}

\begin{abstract}
 Using changes of probability measure  developed by \mbox{Grama} and Haeusler (Stochastic Process.\ Appl., 2000), we obtain two generalizations of the deviation inequalities of Lanzinger and  Stadtm\"{u}ller (Stochastic Process.\ Appl., 2000) and Fuk and Nagaev (Theory Probab. Appl., 1971) to the case of  martingales. Our inequalities recover the best possible decaying rate of independent case.
Applications to linear regressions and weak invariance principles  for martingales are  provided.
\end{abstract}

\begin{keyword} Large deviations; Martingales; Deviation inequalities; Linear regressions; Weak invariance principles
\vspace{0.3cm}
\MSC Primary 60G42, 60E15, 60F10;  Secondary   62J05, 60F17.
\end{keyword}

\end{frontmatter}

%%
%% Start line numbering here if you want
%%
% \linenumbers

%% main text
\section{Introduction}
Assume that $(\xi_i)_{i\geq 1}$ is a sequence of independent and identically distributed (i.i.d.)\ random variables
satisfying the following subexponential condition: for a constant  $\alpha \in (0, 1),$
\begin{eqnarray} \label{B1.3}
 K := \mathbf{E}[ \xi_1^2 \exp\{ (\xi_1^+)^\alpha \} ]  < \infty,
\end{eqnarray}
where $x^+= \max\{x, 0\}.$  Denote by $S_n=\sum_{i=1}^{n}\xi_i$ the partial sums of $(\xi_i)_{i\geq 1}.$
 Lanzinger and  Stadtm\"{u}ller \cite{LS00} have obtained the following subexponential inequality:  for any $x, y>0,$
\begin{eqnarray}\label{B1.4}
\mathbf{P}\Big(S_n \geq  x \Big)   & \leq &  \exp\left\{- \frac{x}{y^{1-\alpha}}\Big( 1- \frac{n K}{2xy^{1-\alpha}} \Big)  \right\}+ \frac{n}{\displaystyle e^{ y^{\alpha}}} \mathbf{E} [\exp\{ (\xi_1^+)^\alpha \} ].
\end{eqnarray}
In particular, by taking $y=x,$ inequality (\ref{B1.4}) implies that  for any $x >0,$
\begin{eqnarray}\label{sfsd}
 \limsup_{n\rightarrow \infty} \frac{1}{n^\alpha} \log \mathbf{P}\Big(S_n \geq  nx \Big)   & \leq & - x^\alpha
\end{eqnarray}
and
\begin{eqnarray}\label{Bsasaf}
 \mathbf{P}\Big(S_n \geq  n  \Big)   & = & O\Big( \exp\Big\{ - c \, n^\alpha \Big\} \Big),   \  \  \ n \rightarrow \infty,
\end{eqnarray}
where $c > 0$ does not depend on $n.$  The last two results (\ref{sfsd}) and (\ref{Bsasaf}) are the best possible under the present condition, since a large deviation principle (LDP)  with good rate function $x^\alpha$  can be obtained in situations where some more information on the tail behavior of $\xi_1$ is available; see Nagaev \cite{N69}.  Under the subexponential condition, more precise estimations on tail probabilities, or large deviation expansions,  can be found in Nagaev \cite{N69,N79}, Saulis and Statulevi\v{c}ius \cite{SS78} and  Borovkov \cite{B00a,B00b}.

%Recently, the generalizations of (\ref{sfsd}) and (\ref{Bsasaf}) have attracted certain interest.
 %When the  coefficients are of subexponential mixing rates, Merlev\`{e}de, Peligrad and Rio \cite{R10} have obtained a  %subexponential inequality of the type (\ref{Bsasaf}) for sums of weakly dependent sequences with summands having subexponential tails. For weighted sum of  i.i.d.\ random variables, we refer to  Gantert, Ramanan and Rembart \cite{GRF14}, where they gave a generalization of (\ref{sfsd}).

%Doukhan  and Neumann \cite{PN07} have established a generalization of (\ref{Bsasaf})
%by introducing a new concept of weak dependence which is more general than mixing. Such conditions are particularly well suited %for deriving estimates for the cumulants of sums of random variables.

Our first aim is to give a generalization of (\ref{Bsasaf}) for   martingales. Let $(\xi_i, \mathcal{F}_{i})_{i\geq 1}$ be a sequence of  martingale differences. Under the Cram\'{e}r condition $\sup_{ i }\mathbf{E}[ \exp\{|\xi_{i}|\}] < \infty$,
Lesigne and Voln\'{y} \cite{LV01} firstly proved that (\ref{Bsasaf}) holds with $\alpha=1/3,$  and that the power $1/3$ is optimal even for the class of stationary and ergodic sequence of martingale differences. Later,  Fan, Grama and Liu \cite{Fx1} generalized the result of Lesigne and Voln\'{y} by proving that (\ref{Bsasaf}) holds under the moment condition $\sup_{ i }\mathbf{E}[\exp\{ |\xi_{i}|^{\frac{2\alpha}{1-\alpha}}\} ]< \infty,$ and that the power $\alpha$ in (\ref{Bsasaf}) is optimal for the class of stationary sequence of martingale differences. It is obvious that the condition $\sup_{ i }\mathbf{E}[\exp\{ |\xi_{i}|^{\frac{2\alpha}{1-\alpha}}\} ]< \infty$ is much stronger than condition (\ref{B1.3}).
Thus, the  result of Fan, Grama and Liu \cite{Fx1} does not imply (\ref{Bsasaf}) in the independent case.

To fill this gap,  we consider the case of the  martingale  differences having bounded conditional subexponential moments. Under this assumption, we can recover the inequalities (\ref{B1.4}), (\ref{sfsd}) and (\ref{Bsasaf}); see  Theorem \ref{th21}. Our first result implies that if $$u:=\max\bigg\{\Big|\Big|\sum_{i=1}^n \mathbf{E} [\xi_i^2 \exp\{ (\xi_{i}^+)^{\alpha} \} |\mathcal{F}_{i-1}]\Big|\Big|_{\infty}, \ \ \  1 \bigg\}  < \infty,$$
then we have for any $x>0,$
\begin{eqnarray}\label{sfdf5}
  \mathbf{P}\Big( \max_{1\leq k \leq n} S_k \geq x   \Big)
 &\leq &  2 \exp\Bigg\{ -  \frac{x^2}{2(u +x^{2-\alpha}) }   \Bigg\}.
\end{eqnarray}
To illustrate our result, consider the simple case that $(\xi_i)_{i\geq 1}$ are i.i.d..\  Then we have $u=O(n)$ as $n\rightarrow \infty.$ It is interesting to see that when $ 0\leq x = o(n^{1/(2-p)}),$ our bound (\ref{sfdf5}) is sub-Gaussian  $\exp\{-x^2/(2 u)\}$, and then it is very tight.   When  $x=ny$ with $y>0$ fixed, our bound (\ref{sfdf5}) is subexponential $ \exp\big\{ - c_y n^{\alpha}  \big\},$ where $c_y >0$ does not depend on $n$. This coincides with (\ref{Bsasaf}).  Moreover, we find that even if $(\xi_i)_{i\geq 1}$  is not stationary, more precisely $u=o(n^{2-\alpha})$, inequality (\ref{sfsd}) is also true; see Corollary  \ref{th1}.

For the methods, an  approach for obtaining subexponential bound  is to combine the method of Lanzinger and  Stadtm\"{u}ller \cite{LS00} and the tower property of conditional expectation.  This approach has been applied in  Fuk \cite{Fu73}, Liu and  Watbled \cite{Liu09a} and Dedecker  and Fan \cite{DF15}. With this approach, one can  obtain inequality (\ref{B1.4}) with  $$nK =  \sum_{i=1}^n \Big|\Big| \mathbf{E} [\xi_i^2 \exp\{ (\xi_{i}^+)^{\alpha} \} |\mathcal{F}_{i-1}]\Big|\Big|_{\infty}.$$
 However, this result is not the best possible in some cases.
In this paper, we introduce a better method based on changes of probability measure  developed by \mbox{Grama} and Haeusler \cite{GH00} (see also \cite{FGL12}). With this method, we  obtain inequality (\ref{B1.4}) with  $$nK =  \Big|\Big| \sum_{i=1}^n\mathbf{E} [\xi_i^2 \exp\{ (\xi_{i}^+)^{\alpha} \} |\mathcal{F}_{i-1}]\Big|\Big|_{\infty}.$$
Since  the later $nK$ is less than the former one, i.e.,
$$   \Big|\Big| \sum_{i=1}^n \mathbf{E} [\xi_i^2 \exp\{ (\xi_{i}^+)^{\alpha} \} |\mathcal{F}_{i-1}]\Big|\Big|_{\infty} \ \leq \ \sum_{i=1}^n \Big|\Big| \mathbf{E} [\xi_i^2 \exp\{ (\xi_{i}^+)^{\alpha} \} |\mathcal{F}_{i-1}]\Big|\Big|_{\infty},$$ our method has certain significant advantage.

As first example to illustrate this advantage, consider the case of self-normalized deviations.   Assume that $(\varepsilon_{i})_{i=1,...,n}$ is a sequence of independent, unbounded  and symmetric around $0$ random variables. Denote by $\xi_i  =  \varepsilon_{i} \Big/\sqrt{\sum_{i=1}^n \varepsilon_{i} ^2}   $
 and $\mathcal{F}_i =\sigma\{ \varepsilon_{j},   1\leq j \leq i,   |\varepsilon_{k}|, 1\leq k \leq n   \}$. Then $(\xi_{i} , \mathcal{F}_i )_{i=1,...,n}$ is also a sequence of martingale differences.   It is easy to see that
 $$    \Big|\Big|\sum_{i=1}^n \mathbf{E} [  (\xi_i ) ^2 \exp\{ ((\xi_i )^+)^{\alpha} \} |\mathcal{F} _{i-1}] \Big|\Big|_{\infty} \, \leq \,  \Big|\Big|\sum_{i=1}^n\frac{ \varepsilon_{i}^2}{\sum_{i=1}^n \varepsilon_{i}^2 }   e \Big|\Big|_{\infty} \, = \,  e,$$
and that, by the fact that $(\varepsilon_{i})_{i=1,...,n}$ are unbounded,
\begin{eqnarray*}
\sum_{i=1}^n \Big|\Big|\mathbf{E} [  (\xi_i ) ^2 \exp\{ ((\xi_i )^+)^{\alpha} \} |\mathcal{F} _{i-1}] \Big|\Big|_{\infty}  \, \geq \, \sum_{i=1}^n \bigg|\bigg|\frac{ \varepsilon_{i}^2}{\sum_{i=1}^n \varepsilon_{i}^2 }  \bigg|\bigg|_{\infty}   \, = \,   n.
\end{eqnarray*} 

Second example  illustrates this advantage.    Assume that $(\varepsilon_{i})_{i=1,...,n}$ is a sequence of independent and unbounded random variables, and that $(\varepsilon_{i})_{i=1,...,n}$ is independent of $(\xi_i, \mathcal{F}_{i})_{i=1,...,n}.$  Assume that $$ \Big| \Big|\mathbf{E} [\xi_i^2 |\mathcal{F}_{i-1} ] \Big| \Big|_{\infty} \geq 1 \ \ \ \  \textrm{and} \ \ \ \  \Big|\Big|\mathbf{E} [\xi_i^2 \exp\{ |\xi_{i}|^{\alpha} \} | \mathcal{F}_{i-1} ] \Big|\Big|_{\infty} \leq D$$ for a constant $D$ and all $i=1,...,n$. Denote by $\xi_i' = \xi_{i}\varepsilon_{i} \Big/\sqrt{\sum_{i=1}^n \varepsilon_{i} ^2}   $
 and $\mathcal{F}_i'=\sigma\{ \varepsilon_{j},   1\leq j \leq i, \mathcal{F}_{i}, |\varepsilon_{k}|, 1\leq k \leq n   \}$. Then $(\xi_{i}', \mathcal{F}_i')_{i=1,...,n}$ is also a sequence of martingale differences.   It is easy to see that
 $$    \Big|\Big|\sum_{i=1}^n \mathbf{E} [  (\xi_i') ^2 \exp\{ ((\xi_i')^+)^{\alpha} \} |\mathcal{F}'_{i-1}] \Big|\Big|_{\infty} \, \leq \,  \Big|\Big|\sum_{i=1}^n\frac{ \varepsilon_{i}^2}{\sum_{i=1}^n \varepsilon_{i}^2 }    \mathbf{E} [\xi_i^2 \exp\{ |\xi_i|^{\alpha} \} |\mathcal{F}_{i-1} ]\Big|\Big|_{\infty} \, \leq \,  D,$$
and that, by the fact that $(\varepsilon_{i})_{i=1,...,n}$ are unbounded,
\begin{eqnarray*}
\sum_{i=1}^n \Big|\Big|\mathbf{E} [  (\xi_i') ^2 \exp\{ ((\xi_i')^+)^{\alpha} \} |\mathcal{F}'_{i-1}] \Big|\Big|_{\infty}  \, \geq \, \sum_{i=1}^n \bigg|\bigg|\frac{ \varepsilon_{i}^2}{\sum_{i=1}^n \varepsilon_{i}^2 }  \mathbf{E} [\xi_i^2 |\mathcal{F}_{i-1} ]\bigg|\bigg|_{\infty} \, \geq \, \sum_{i=1}^n \bigg|\bigg|\frac{\varepsilon_{i}^2 }{ \sum_{i=1}^n \varepsilon_{i} ^2}\bigg|\bigg|_{\infty}  \, = \,   n.
\end{eqnarray*}
Thus the advantage of our method  is significant.

With changes of probability measure, we also generalize  the following inequality  of Fuk  for martingales (cf. Corollary $3'$ of Fuk \cite{Fu73}; see also Nagaev \cite{N79} for independent case): if ${\mathbf E} [|\xi_i|^p|\mathcal{F}_{i-1}]   < \infty $ for a $p\geq 2$ and all $i \in [1, n]$, then for any $x  > 0$,
\begin{eqnarray}  \label{s106}
 {\mathbf P}\bigg( \max_{1\leq k \leq n} S_k \geq x   \bigg)   \   \leq \  \exp\Bigg\{- \frac{ x^2}{2 \widetilde{V}^2 } \Bigg \} +  \frac{ \widetilde{C}_p}{x^p}   \,  ,
\end{eqnarray}
where
$$\widetilde{V}^2  := \frac{1}{4}(p+2)^2 e^p\, \sum_{i=1}^{n}  \Big|\Big|\mathbf{E} [\xi_i^2 |\mathcal{F}_{i-1} ] \Big|\Big|_{\infty}\ \   \  \textrm{and}  \ \ \   \widetilde{C}_p:=  \Big( 1+ \frac{2}{p} \Big)^p \, \sum_{i=1}^{n}\Big|\Big| \mathbf{E} [|\xi_i|^p |\mathcal{F}_{i-1} ] \Big|\Big|_{\infty}.$$
In Corollary \ref{co22}, we prove that (\ref{s106}) holds true when $\widetilde{V}^2$ and $\widetilde{C}_p$ are replaced by the following two smaller values $V^2$ and $C_p$  respectively, where
$$V^2  := \frac{1}{4}(p+2)^2 e^p\, \Big|\Big|\sum_{i=1}^{n}  \mathbf{E} [\xi_i^2 |\mathcal{F}_{i-1} ] \Big|\Big|_{\infty}\ \   \  \textrm{and}  \ \ \   C_p:=  \Big( 1+ \frac{2}{p} \Big)^p \, \Big|\Big|\sum_{i=1}^{n} \mathbf{E} [|\xi_i|^p |\mathcal{F}_{i-1} ] \Big|\Big|_{\infty}.$$
To illustrate the improvement of Corollary \ref{co22} on Fuk's inequality (\ref{s106}), consider the following comparison between $C_p$ and $\widetilde{C}_p $ in the case of random weighted self-normalized deviations. As before, assume that $(\varepsilon_{i})_{i=1,...,n}$ is a sequence of independent and unbounded random variables, and that $(\varepsilon_{i})_{i=1,...,n}$ is independent of $(\xi_i, \mathcal{F}_{i})_{i=1,...,n}.$ Assume   $$ 1\leq  \Big|\Big| \mathbf{E} [ |\xi_i|^p   | \mathcal{F}_{i-1} ]\Big|\Big|_{\infty}  \leq E$$ for a constant $E$ and all $i=1,...,n$.   Denote by $\xi_i'  = \xi_{i}\varepsilon_{i} /(\sum_{i=1}^n |\varepsilon_{i} |^p )^{1/p} $
 and $\mathcal{F}_i'=\sigma\{  \mathcal{F}_{i}, |\varepsilon_{k}|, 1\leq k \leq n   \}$. Then $(\xi_{i}', \mathcal{F}_i')_{i=1,...,n}$ is also a sequence of martingale differences.   It is easy to see that
 $$  \Big|\Big|  \sum_{i=1}^{n}  \mathbf{E} [|\xi_i'|^p |\mathcal{F}_{i-1}' ] \Big|\Big|_{\infty}\ \leq \ \Big|\Big|\sum_{i=1}^{n} \frac{|\varepsilon_{i}|^p }{ \sum_{i=1}^n |\varepsilon_{i} |^p}   \mathbf{E} [|\xi_i |^p |\mathcal{F}_{i-1} ] \Big|\Big|_{\infty}\ \leq  \   E,$$
and that, by the fact that $(\varepsilon_{i})_{i=1,...,n}$ are unbounded,
\begin{eqnarray*}
\sum_{i=1}^n \Big|\Big| \mathbf{E} [|\xi_i'|^p |\mathcal{F}_{i-1}' ] \Big|\Big|_{\infty}  \, \geq \, \sum_{i=1}^n \bigg|\bigg|\frac{|\varepsilon_{i}|^p }{ \sum_{i=1}^n |\varepsilon_{i} |^p}   \mathbf{E} [|\xi_i |^p |\mathcal{F}_{i-1} ] ]\bigg|\bigg|_{\infty}   \, \geq \, n.
\end{eqnarray*}
Hence $C_p$ is much smaller than  $\widetilde{C}_p.$ The improvement of Corollary \ref{co22} on Fuk's inequality (\ref{s106}) is significant.

For two positive constants $\delta$ and $C,$ assume either $ \mathbf{E} [|\xi_i|^{p+\delta} ] \leq C$ and $ \mathbf{E} [\xi_i^2 |\mathcal{F}_{i-1} ] \leq C $ a.s.\ or  $ \mathbf{E} [|\xi_i|^p |\mathcal{F}_{i-1} ] \leq C$ a.s.\ for a $p\geq 2$ and all $i=1,...,n.$
Then we have  for  any $\alpha \in (\frac12, \infty),$
\begin{eqnarray} \label{lsvss}
{\mathbf P}\Big( \max_{1\leq k \leq n} S_k \geq n^\alpha     \Big)   \   = \   O \Big(\frac{1 }{n^{\alpha p-1}} \Big) \, ,\ \ \ \ \ \  \ n \rightarrow \infty,
\end{eqnarray}
See Theorem  \ref{th23} and Corollary \ref{co22}. Under a stronger condition that $(\xi_i)_{i=1,...,n}$ have bounded conditional moments, inequality (\ref{lsvss}) improves a result of  Lesigne and Voln\'{y} \cite{LV01} where   Lesigne and Voln\'{y}
proved that if $ \mathbf{E} [|\xi_i|^p   ] \leq C $ for a constant $C,$ then
\begin{eqnarray}
 {\mathbf P}\Big(  S_n \geq n    \Big)   \   = \  O\Big(  \frac{1 }{n^{ p/2}} \Big)\, ,\ \ \ \ \ \ \ \ \ \  \ n \rightarrow \infty,
\end{eqnarray}
 and that the order $n^{ -p/2}$ of the last inequality   is optimal even for the class of stationary and ergodic sequence of martingale differences.

The paper is organized as follows. We present our main results in Section \ref{sec2}, and discuss the applications  to linear regressions and  weak invariance principle  in Section \ref{sec3}. The proofs of theorems are given in Sections  \ref{sec4} - \ref{sec7}.

\section{Main results}\label{sec2}

Assume that we are given a sequence of real-value martingale differences $(\xi _i,\mathcal{F}_i)_{i=0,...,n}$
 defined on some probability space $(\Omega ,\mathcal{F},\mathbf{P})$, where $\xi _0=0 $ and $\{\emptyset, \Omega\}=%
\mathcal{F}_0\subseteq ...\subseteq \mathcal{F}_n\subseteq \mathcal{F}$ are
increasing $\sigma$-fields. So we have $\mathbf{E}[\xi_{i}|\mathcal{F}_{i-1}]= 0, \  i=1,...,n, $ by definition. Set
\begin{equation}  \label{matingal}
S_0=0 \ \ \ \ \ \ \textrm{and} \ \ \ \ \ S_k=\sum_{i=1}^k\xi _i,\ \  \ \  \quad k=1,...,n.
\end{equation}
Then $S:=(S_k,\mathcal{F}_k)_{k=0,...,n}$ is a martingale. Let $\left\langle S\right\rangle $ be the quadratic characteristic of the
martingale $S:$
\begin{equation}\label{quad}
\left\langle S\right\rangle _0=0 \ \ \ \ \  \ \textrm{and}\ \ \  \ \ \    \left\langle S\right\rangle _k=\sum_{i=1}^k\mathbf{E}[\xi _i^2|\mathcal{F}
_{i-1}],\quad k=1,...,n.
\end{equation}

Our first result is the following subexponential inequality on tail probabilities for martingales. A similar inequality for separately Lipschitz functionals has been obtained recently by Dedecker  and Fan \cite{DF15}.
\begin{theorem}\label{th21}
Assume $$  C_n :=\sum_{i=1}^n \mathbf{E}[ \xi_i^2 \exp\{ (\xi_{i}^+)^{\alpha} \}]< \infty$$ for a  constant $\alpha \in (0,1)$. Denote by
\begin{eqnarray}\label{sads}
\Upsilon(S)_k= \sum_{i=1}^k \mathbf{E} [\xi_i^2 \exp\{ (\xi_{i}^+)^{\alpha} \} |\mathcal{F}_{i-1}],\ \ \ \ \ \ \  \ k \in [1,n].
\end{eqnarray}
Then  for all $x, u > 0$,
\begin{eqnarray}\label{ineq5}
 &&  \mathbf{P}\Big(  S_k \geq x \ \mbox{and}\ \Upsilon(S)_k \leq u \ \mbox{for some}\ k \in [1,n] \Big) \ \ \ \ \  \ \ \  \nonumber \\
&& \ \ \leq \left\{ \begin{array}{ll}
\exp\Bigg\{\displaystyle -  \frac{x^2}{2u }   \Bigg\} +  C_n \bigg( \frac{ x  } {u }\bigg)^{2/(1-\alpha)}\exp\Bigg\{ -\Big(\frac{u}{x  } \Big)^{\alpha/(1-\alpha)} \Bigg\} \ \  & \textrm{if\ \  $0\leq x < u^{1/(2-\alpha)}$}   \\
\vspace{-0.2cm}\\
\exp \Bigg\{\displaystyle - x^\alpha \Big( 1- \frac{u}{2\, x^{2-\alpha}}\Big)  \Bigg\} + C_n \frac{1}{x^2}\exp \Bigg\{- x^{\alpha}   \Bigg\}   & \textrm{if\ \  $x \geq u^{1/(2-\alpha)}$. }
\end{array} \right.
\end{eqnarray}
\end{theorem}

It is obvious  that $$C_n = \mathbf{E}[ \Upsilon(S)_n ] \leq ||\Upsilon(S)_n||_{\infty}.$$ Hence, if $u\geq \max\{||\Upsilon(S)_n||_{\infty}, 1 \},$  then (\ref{ineq5}) implies the following rough bounds
\begin{eqnarray}
  \mathbf{P}\Big( \max_{1\leq k \leq n} S_k \geq x   \Big)
  &\leq& \left\{ \begin{array}{ll}
2 \exp\bigg\{\displaystyle -  \frac{\ x^2}{2\, u }   \bigg\} \ \ \ & \textrm{if\ \ $0\leq x < u^{1/(2-\alpha)}$}  \\
\vspace{-0.3cm}  \\
2 \exp \bigg\{\displaystyle  - \frac12 x^\alpha   \bigg\} \ \ \ & \textrm{if\ \ $x \geq u^{1/(2-\alpha)}$ }
\end{array} \right. \label{ineq6}\\
 &\leq &  2 \exp\Bigg\{ -  \frac{x^2}{2(u +x^{2-\alpha}) }   \Bigg\}. \label{ineq7}
\end{eqnarray}
Thus for moderate $ x \in (0,  u^{1/(2-\alpha)})$,
 bound (\ref{ineq5}) is sub-Gaussian.
For  all $x\geq u^{1/(2-\alpha)},$ bound (\ref{ineq5}) is subexponential, and is of the order $ \exp\Big\{ - \frac{1}2 x^{\alpha} \Big\}.$
Moreover, when $\frac{x}{u^{1/(2-\alpha)} }   \rightarrow \infty,$ by (\ref{ineq5}), this order can be improved to $\exp\Big\{ - (1+\varepsilon)x^{\alpha}  \Big\}$ for any given $\varepsilon>0$.

Define
\begin{eqnarray*}
\widehat{\Upsilon} (S)_k= \sum_{i=1}^k \mathbf{E} [\xi_i^2 \exp\{ |\xi_{i}|^{\alpha} \} |\mathcal{F}_{i-1}]
\end{eqnarray*} for any $k \in [1,n]$. Then it holds $ \Upsilon  (S)_k \leq \widehat{\Upsilon} (S)_k.$
It is obvious that  bound (\ref{ineq5}) is also the upper bound on the tail probabilities
$$\mathbf{P} \Big( \pm S_k \geq x \ \mbox{and}\ \widehat{\Upsilon} (S)_k \leq u \ \mbox{for some}\ k \in [1,n] \Big).$$ Moreover, if $||\widehat{\Upsilon} (S)_n||_{\infty} \leq u,$ then bound (\ref{ineq5}) is an upper bound on the partial sums  tail probabilities
$\mathbf{P}( \pm \max_{1\leq k \leq n} S_k \geq x ).$

When $(\xi_i)_{i\geq 1}$ are i.i.d.\ random variables, then we have $C_n =\Upsilon(S)_n=O(n)$ as $n\rightarrow \infty$. In this case, inequality (\ref{ineq7}) implies the following  large deviation inequality: for any $x>0,$
\begin{eqnarray}\label{fineq6}
  \mathbf{P}\Big( \max_{1\leq k \leq n} S_k \geq n x   \Big)  \ = \ O \Big( \exp \Big\{ - c_x n^\alpha  \Big\} \Big), \ \ \ \ n\rightarrow \infty,
\end{eqnarray}
where $c_x> 0$ does not depend  on $n$. Moreover,  the following  LDP   result for
martingales shows that $c_x$ in (\ref{fineq6}) is close to $x^\alpha$, and that  $||\Upsilon(S)_n||_{\infty}$ is allowed to tend to infinity in an order larger than $n$.
\begin{corollary}\label{th1}
Assume condition  of Theorem \ref{th21}.
If
\begin{eqnarray}
||\Upsilon(S)_n||_{\infty} \ = \ o(n^{2-\alpha}),\ \ \  \ n \rightarrow \infty,
\end{eqnarray}
then for any $x\geq 0$,
\begin{eqnarray}
 \limsup_{n\rightarrow \infty}\frac{1}{n^\alpha}\log \mathbf{P}\left( \max_{1\leq k \leq n} \frac{1}{n}S_k \geq x \right) \ \leq \  -  x ^\alpha . \,
\end{eqnarray}
\end{corollary}

\begin{remark}  \label{re1}
This result cannot be  improved under the present condition even for the class of i.i.d.\ random variables; see Nagaev \cite{N69}. In fact, if $(\xi_i)_{i\geq 1}$ are i.i.d.\ and satisfy the following condition  for an integer $p \geq 2$ and all $x$ large enough,
\begin{eqnarray}\label{as11}
 \frac{1}{ x^{2p}} \, \exp \Big\{ - x^\alpha  \Big\} \ \leq \  \mathbf{P}\Big( |\xi_1| \geq   x   \Big)\  \leq\  \frac{1}{ x^{3}}  \, \exp \Big\{ - x^\alpha  \Big\},
\end{eqnarray}
then we have
$\Upsilon(S)_n   = o(n^{2-\alpha})$ as $n\rightarrow \infty$ and for any $x>0,$
\begin{eqnarray}\label{fineq15}
 \lim_{n\rightarrow \infty}\frac{1}{n^\alpha}\log \mathbf{P}\left( \max_{1\leq k \leq n} \frac{1}{n}S_k \geq x \right)\  = \ -  x ^\alpha.
\end{eqnarray}
\end{remark}

If the martingale differences $(\xi _i,\mathcal{F}_i)_{i=0,...,n}$ have $p$-th moments ($p\geq 2$), then we have the following
inequality, which is  similar to the results of  Haeusler \cite{H84} and   \cite{FGL12}.

\begin{theorem}\label{th22}
 Let $p\geq 2$. Assume ${\mathbf E} [|\xi_i|^p]   < \infty $ for all $i \in [1, n]$.  Denote by
$$ \Xi(S)_k \ = \ \sum_{i=1}^{k} \mathbf{E} [(\xi_i^+)^p |\mathcal{F}_{i-1} ],\ \ \ \ \ \ \ \ k \in [1, n]. $$
Then  for all $x, y, v, w  > 0$,
\begin{eqnarray}\label{sdfds}
&& {\mathbf P}\Big( S_k \geq x,\ \langle S\rangle_k \leq v  \  \textrm{and} \ \ \Xi(S)_k   \leq  w \textrm{ for some }k \in [1, n] \Big) \nonumber  \\ &&\ \ \ \ \ \   \leq \  \exp\Bigg\{\displaystyle - \frac{\alpha^2x^2}{2e^p v} \Bigg \} + \exp\Bigg\{\displaystyle - \frac{\beta x}{y} \log\Big( 1 + \frac{\beta x y^{p-1}}{w}\Big)  \Bigg \}+ \mathbf{P}\left( \max_{ 1\leq i\leq n}\xi_{i} > y \right)  ,\ \ \ \ \ \  \
\end{eqnarray}
where
\begin{eqnarray}\label{defab}
\alpha= \frac{2}{p+2}\ \ \ \ \  \ \textrm{and} \ \ \  \ \ \ \beta=1-\alpha\, .
\end{eqnarray}
\end{theorem}

Setting $y=\beta x,$ we obtain the following generalization of the Fuk-Nagaev  inequality (\ref{s106}).
\begin{corollary}\label{co22}
 Let $p\geq 2$. Assume $\big|\big|{\mathbf E} [|\xi_i|^p |\mathcal{F}_{i-1}] \big|\big|_{\infty}  < \infty $ for all $i \in [1, n]$.
It holds for all $x  > 0$,
\begin{eqnarray}\label{sfddghfd}
 {\mathbf P}\Big( \max_{1\leq k \leq n} S_k \geq x   \Big)   \   \leq \  \exp\Bigg\{- \frac{ x^2}{2\, V^2 } \Bigg \} +  \frac{ C_p}{x^p}   \, ,\ \ \ \ \ \  \
\end{eqnarray}
where
\begin{eqnarray}\label{cpdef}
V^2=\frac{1}{4}(p+2)^2 e^p\, \Big|\Big|\left\langle S\right\rangle_n\Big|\Big|_{\infty} \ \ \ \ \  \textrm{and}\ \ \ \ \ \  C_p = \Big( 1+ \frac{2}{p} \Big)^p \, \bigg|\bigg| \sum_{i=1}^{n} \mathbf{E} [|\xi_i|^p |\mathcal{F}_{i-1} ] \bigg|\bigg|_{\infty}.
\end{eqnarray}
\end{corollary}

It is worth noting that if $ \mathbf{E} [|\xi_i|^p |\mathcal{F}_{i-1} ] \leq C$ for a constant $C$ and  all $i \in [1, n]$,   then, by Jensen's inequality, it holds  $ \mathbf{E} [\xi_i^2 |\mathcal{F}_{i-1} ] \leq C^{2/p}$ for all $i \in [1, n].$ Inequality (\ref{sfddghfd}) implies  the following sub-Gaussian bound  for  any  $ x = O(\sqrt{n} \, (\ln n)^\beta ), $ $n\rightarrow \infty,$ with $ \beta$ satisfying $ \beta>0$ if $p=2$ and $ \beta \in (0, 1/2]$  if $p > 2,$
\begin{eqnarray}\label{sdfg03}
 {\mathbf P}\Big( \max_{1\leq k \leq n} S_k \geq   x \Big)   \   =  \   O\bigg( \exp\bigg\{- C \, \frac{ x^2}{n} \bigg \} \bigg),
\end{eqnarray}
where $C>0$ does not depend on $x$ and $n.$
The bound (\ref{sdfg03}) is similar to the classical Azuma-Hoeffding inequality, and thus it is tight.
Inequality (\ref{sfddghfd})  also implies that
for any $\alpha \in (\frac12, \infty)$ and any $x>0,$
\begin{eqnarray}\label{sdfg}
 {\mathbf P}\Big( \max_{1\leq k \leq n} S_k > n^\alpha x   \Big)     &=&    O\Big(  \frac{c_x}{n^{\alpha p-1}} \Big) \, ,\ \ \ \ \ \  \  n\rightarrow \infty,
\end{eqnarray}
where $c_x > 0$ does not depend on $n.$ Equality (\ref{sdfg}) is first obtained by   Fuk   \cite{Fu73} and it is the best possible under the stated condition even for the sums of independent  random variables (cf.\ Fuk and Nagaev \cite{FN71}).

If the martingale differences $(\xi _i,\mathcal{F}_i)_{i=1,...,n}$ satisfy $ \mathbf{E} [|\xi_i|^p   ] \leq C $ for a constant $C$ and all $i \in [1, n],$ then Lesigne and Voln\'{y} \cite{LV01} proved  that
for  any $x>0,$
\begin{eqnarray}\label{lvs}
 {\mathbf P}\Big(   S_n > n  x   \Big)  \  =  \  O\Big(   \frac{c_x}{n^{ p/2}}  \Big)\, ,\ \ \ \ \ \  \  n\rightarrow \infty,
\end{eqnarray}
where $c_x>0$ does not depend on $n,$  and that the order $n^{ -p/2}$ is optimal even for the class of stationary and ergodic sequence of martingale differences. When $\alpha =1,$ equality (\ref{sdfg}) implies  the following large deviation convergence rate  for  any $x>0,$
\begin{equation}\label{sfdf01}
 {\mathbf P}\Big( \max_{1\leq k \leq n} S_k > n  x   \Big)   \  =  \  O\Big(  \frac{c_x}{n^{p-1}} \Big)\, ,\ \ \ \ \ \  \  n\rightarrow \infty,
\end{equation}
where $c_x>0$ does not depend on $n.$ When $p\geq2,$ it holds $p-1 \geq p/2.$
Thus (\ref{sfdf01}) refines the bound (\ref{lvs})   under the stronger assumption that the $p\,$-th
conditional moments are uniformly bounded. Moreover, the following proposition of  Lesigne and Voln\'{y} \cite{LV01} shows that  the estimate of (\ref{sfdf01}) cannot be essentially improved even in the i.i.d.\ case.
\vspace{0.2cm}

\noindent\textbf{Proposition A.} \emph{Let $p\geq 1$ and $(c_n)_{n\geq 1}$ be a real positive sequence approaching zero.
There exists a sequence of  i.i.d.\ random variables $(\xi _i)_{i\geq 1}$ such that $\mathbf{E}[|\xi_i|^p] < \infty,$ $\mathbf{E}[\xi_i] =0$ and}
\[
\limsup_{n \rightarrow \infty} \frac{n^{p-1}}{c_n}\mathbf{P}(|S_n| \geq n) = \infty.
\]

When ${\mathbf E} [|\xi_i|^2 |\mathcal{F}_{i-1}]$ and ${\mathbf E} [|\xi_i|^p], $ for a $p> 2$ and all $i=1,...,n,$ are  all uniformly bounded  (but the condition $\mathbf{E} [|\xi_i|^p |\mathcal{F}_{i-1} ] \leq C$ may be violated  for some $i \in [1, n]$),   we have the following result.
\begin{theorem} \label{th23}Let $p\geq 2$. Assume ${\mathbf E} [|\xi_i|^{p+\delta}]   < \infty $ for a small $\delta>0$ and all $i \in [1, n]$.
Then  for all $x,   v > 0$,
\begin{eqnarray}\label{ineq1ss3}
&&{\mathbf P}\Big( S_k \geq x  \  \textrm{and} \ \langle S\rangle_k \leq v^2   \textrm{ for some }k \in [1, n] \Big) \ \ \ \ \ \ \ \ \ \nonumber \\
 &&\ \ \ \ \ \ \ \ \ \  \ \ \leq \ \exp\left\{  - \frac{x^2}{2 \Big(v^2+ \frac{1}{3}x^{(2p+\delta)/(p+\delta) } \Big)} \right\}+ \frac{1}{x^{p}} \sum_{i=1}^n{\mathbf E} \Big[|\xi_i|^{p+\delta}\mathbf{1}_{\{\xi_{i}>x^{p/(p+\delta) }\}} \Big].
\end{eqnarray}
\end{theorem}

If $ \mathbf{E} [|\xi_i|^{p+\delta} ] \leq C$ and $ \mathbf{E} [\xi_i^2 |\mathcal{F}_{i-1} ] \leq C $ for a constant $C$ and all $i \in [1,n]$, then (\ref{ineq1ss3}) implies that  for  any $\alpha \in (\frac12, \infty)$ and any $x>0,$
\begin{eqnarray}\label{fsineq1}
 {\mathbf P}\Big( \max_{1\leq k \leq n} S_k \geq n^\alpha x   \Big)     &\leq &  \exp\bigg\{- c_x \, \min\Big\{ n^{\frac{\alpha\delta}{p+\delta} }, n^{2\alpha -1 } \Big\} \bigg\}+ \frac{C/x^p}{n^{\alpha p-1} }  \nonumber \\
  &=&  O\Big(  \frac{c_x}{n^{\alpha p-1}} \Big)  \, ,\ \ \ \  \ \ \ n \rightarrow \infty,
\end{eqnarray}
where $c_x > 0$ does not depend on $n.$ Thus  (\ref{sdfg}) and (\ref{fsineq1}) have the same  convergence rate.

The different between the conditions of (\ref{sdfg}) and (\ref{fsineq1}) is that the assumption $ \mathbf{E} [|\xi_i|^p |\mathcal{F}_{i-1} ] \leq C$ has been replaced by the two  assumptions $ \mathbf{E} [|\xi_i|^{p+\delta} ] \leq C$ and $ \mathbf{E} [\xi_i^2 |\mathcal{F}_{i-1} ] \leq C $ for   all $i \in [1,n]$. Notice that the two assumptions $ \mathbf{E} [|\xi_i|^p |\mathcal{F}_{i-1} ] \leq C$ and $ \mathbf{E} [|\xi_i|^{p+\delta} ] \leq C$ are not included in each other. Thus  Corollary \ref{co22} and Theorem  \ref{th23} do not imply each other.

From Theorem \ref{th23},  the following corollary is obvious.
\begin{corollary}\label{co23} Assume the condition of Theorem \ref{th23}.
Then  for all $x, v > 0$,
\begin{eqnarray}
 {\mathbf P}\Big( \max_{1\leq k \leq n} S_k \geq x   \Big)
  &\leq& \exp\left\{  - \frac{x^2}{2 \big(n v^2+ \frac{1}{3}x^{(2p+\delta)/(p+\delta) } \big)} \right\}+ \frac{1}{x^{p}} \sum_{i=1}^n{\mathbf E} \Big[|\xi_i|^{p+\delta}\mathbf{1}_{\{\xi_{i}>x^{p/(p+\delta) }\}} \Big] \nonumber \\
  && \ \  + \,  \frac{1}{ \, v^{p+\delta }}  \mathbf{E}\Big[\Big|\frac{\langle S\rangle_n}{n} \Big|^{(p+\delta)/2}\Big]  . \label{fineq31}
\end{eqnarray}
Moreover, it holds
\begin{eqnarray}\label{fineq32}
\mathbf{E}\Big[\Big|\frac{\langle S\rangle_n}{n} \Big|^{(p+\delta)/2}\Big]
  &\leq& \frac1n  \sum_{i=1}^n{\mathbf E} \big[|\xi_i|^{p+\delta} \big].
\end{eqnarray}
\end{corollary}

Inequality (\ref{fineq32}) implies that if  $\sup_{i}{\mathbf E} [|\xi_i|^{p}]   < \infty$ for a $p\geq2, $ then $ \mathbf{E}\big[\big| \langle S\rangle_n/n \big|^{p/2}\big]$ are uniformly bounded for all $n$.

Compared to Theorem  \ref{th23}, Corollary \ref{co23} is more applicable since it only need the moment of $\langle S\rangle_n$ instead of the uniform bound of $\langle S\rangle_n$.

%We will make use of Corollary \ref{co23} to prove a weak invariance principle for martingales (cf.\ Theorem  \ref{f1jka}).

Assume  ${\mathbf E} [|\xi_i|^{p+\delta}]  \leq C $ for a  $p\geq 2$ and all $i \in [1,n]$ (without any condition on $\langle S\rangle_n$).
Applying (\ref{fineq32}) to (\ref{fineq31}) with $n v^2=\frac{2}{3}x^{(2p+\delta)/(p+\delta) },$ we have for all $x, v > 0$,
\begin{eqnarray}
 {\mathbf P}\Big( \max_{1\leq k \leq n} S_k \geq x   \Big)
  &\leq&  \exp\left\{  - \frac{1}{2} x^{\delta/(p+\delta) }   \right\}+ \frac{n C}{x^{p}}
  \  + \, \Big( \frac{3n}2\Big)^{\frac{p+\delta}{2} }  \frac{C}{ x^{p+\delta/2 }}.
\end{eqnarray}
The last inequality shows that for any  $x  > 0$,
\begin{eqnarray}\label{cvdsds}
 {\mathbf P}\Big( \max_{1\leq k \leq n} S_k \geq n   \Big)
  & =&  O \Big( \frac{1}{ n^{p/2 }} \Big), \ \ \ \ \ \ \ \ \ n\rightarrow \infty.
\end{eqnarray}
Since $\delta >0$ can be any small, equality (\ref{cvdsds}) is closed to the best possible large deviation convergence rate $n^{-(p+\delta)/2 }$ given by  Lesigne and Voln\'{y} \cite{LV01} (cf.\ (\ref{lvs})).

\section{Applications}\label{sec3}

The exponential concentration inequalities for martingales  have many applications.  \mbox{McDiarmid} \cite{M}, Rio \cite{R13a} and Dedecker and Fan \cite{DF15}  applied such type inequalities to estimate the concentration of separately Lipschhitz functions.  Liu  and  Watbled \cite{Liu09a} adopted these inequalities to deduce asymptotic properties of the free energy of directed polymers in a random environment.  We refer to Bercu  and Touati \cite{BT08} and \cite{FGL15} for more interesting applications of the concentration inequalities for martingales.
In the sequel, we discuss how to apply our results to linear regression models and  weak invariance principles.

\subsection{Linear regressions}\label{sec3.1}
Linear regressions can be used to investigate the impact of one variable on the other, or to predict the value of one variable based on the other. For instance, if one wants to see impact of footprint size on height, or predict height according to a certain given value of footprint size. The stochastic linear regression model is given by, for all $k  \in [1,   n],$
\begin{equation}\label{ine29}
X_{k}=\theta \phi_k + \varepsilon_{k}\, ,
\end{equation}
where $(X_k)_{k=1,...,n}, (\phi_k)_{k=1,...,n}$ and $(\varepsilon_{k})_{k=1,...,n}$ are the observations, the regression variables and the driven noises, respectively. We assume that $(\phi_k)_{k=1,...,n}$ is a sequence of independent  random variables, and that  $(\varepsilon_k )_{k=1,...,n}$ is a sequence of  martingale differences with respect to
the natural filtration. Moreover, we suppose that $(\phi_k)_{k=1,...,n}$ and  $(\varepsilon_k)_{k=1,...,n}$ are independent. Our interest is to estimate the unknown parameter $\theta.$ The well-known  least-squares estimator $\theta_n$ is given below
\begin{equation}\label{ine30}
\theta_n = \frac{\sum_{k=1}^n \phi_{k} X_k}{\sum_{k=1}^n \phi_{k}^2}.
\end{equation}
Recently, Bercu and Touati \cite{BT08} have obtained some very precise exponential bounds on the tail probabilities $\mathbf{P}\left( |\theta_n -\theta|  \geq x  \right).$ However, their precise  bounds  depend   on the distribution of
input random variables $(\phi_k)_{k=1,...,n},$ which restricts the applications of these bounds when the distributions of
input random variables are unknown.
When $(\varepsilon_k)_{k=1,...,n}$ are independent normal  random variables  with a common variation $\sigma^2> 0,$ Liptser and Spokoiny  \cite{Ls01} have established the following estimation: for all $x\geq 1,$
\begin{eqnarray}\label{ls21}
 \mathbf{P}\left( \pm \, (\theta_n -\theta)\sqrt{ \Sigma _{k=1}^n \phi_{k}^2} \geq x  \right) \  \leq\    \sqrt{\frac2 \pi} \, \frac{ \sigma }{x}   \, \exp\Bigg\{ - \frac{x^2}{2 \, \sigma^2  }  \Bigg\}.
\end{eqnarray}
When $(\varepsilon_k)_{k=1,...,n}$ are conditionally sub-Gaussian, similar estimation  is allowed to be obtained in Liptser and Spokoiny  \cite{Ls01}.
An interesting feature of bound (\ref{ls21}) is that the  bound  does not depend  on the distribution of
input random variables.
Here, we would like to   generalize  inequality  (\ref{ls21}) to the case that $(\varepsilon_k)_{k=1,...,n}$ are martingale differences  and also non sub-Gaussian.

\begin{theorem} \label{thlin}
Assume for two  constants $\alpha \in (0,1)$ and $D$,
 \[
\mathbf{E}\Big[ \varepsilon_i^2 e^{ |\varepsilon_{i}|^{\alpha}  } \, \Big| \, \sigma\{\varepsilon_j
, j\leq i-1\}  \Big] \leq D
\]
for all  $i \in [1, n].$ Then for any $u\geq \max\{D, 1 \}$ and all $x>0$,
\begin{eqnarray}
 \mathbf{P}\left(  \pm (\theta_n -\theta)\sqrt{ \Sigma _{k=1}^n \phi_{k}^2} \geq x  \right)
  &\leq& \left\{ \begin{array}{ll}
2 \exp\bigg\{\displaystyle -  \frac{\, x^2}{2\, u }   \bigg\} \ \ \ & \textrm{if\ \ $0\leq x < u^{1/(2-\alpha)}$}  \\
\vspace{-0.3cm}\\
2 \exp \bigg\{\displaystyle  - \frac12 x^\alpha   \bigg\} \ \ \ & \textrm{if\ \ $x \geq u^{1/(2-\alpha)}$ }
\end{array} \right.  \label{th5ineq}\\
 &\leq &  2 \exp\Bigg\{\displaystyle  -  \frac{x^2}{2(u +x^{2-\alpha}) }   \Bigg\}.
\end{eqnarray}
In particular, it holds  for any $x>0,$
\begin{eqnarray}\label{sfdf37}
 {\mathbf P}\bigg( \pm (\theta_n -\theta)\sqrt{\Sigma_{k=1}^n \phi_{k}^2} \geq \sqrt{n}\, x   \bigg)   \   = \  O\bigg( \exp \Big\{- c_x n^{\alpha/2} \Big\}   \bigg),\ \ \ \ \ n\rightarrow \infty,
\end{eqnarray}
where $c_x>0$ does not depend on $n.$
\end{theorem}

If $(\varepsilon_k)_{k=1,...,n}$ have the Weibull distributions
 and the conditional variances are uniformly bounded,    then we have the following inequality which has the same exponentially decaying rate of (\ref{sfdf37}).
\begin{theorem} \label{thFA12}
Assume for three constants $\alpha \in (0,1),$ $E$ and $F,$
 \[
\mathbf{E}\Big[ \varepsilon_i^2  \, \Big| \, \sigma\{\varepsilon_j
, j\leq i-1\}  \Big] \leq E \ \ \ \ \
\textrm{and}
\ \ \ \ \  \mathbf{E}\Big[ \exp\{|\varepsilon_i|^{\frac{\alpha}{1-\alpha}} \}   \Big] \leq F
\]
for all $i \in [1, n].$
 Then for  all $x>0,$
\begin{eqnarray}
 \mathbf{P}\left(  \pm (\theta_n -\theta)\sqrt{ \Sigma _{k=1}^n \phi_{k}^2} \geq x  \right)
\ \leq \  \exp\Bigg\{\displaystyle  -  \frac{x^2}{2(E + \frac13 x^{2-\alpha}) }   \Bigg\} + n F \, \exp\Big\{- x^\alpha \Big\}.
\end{eqnarray}
In particular, equality (\ref{sfdf37})  holds.
\end{theorem}

If $(\varepsilon_k)_{k=1,...,n}$ have finite conditional moments, by Corollary \ref{co22}, then
 we have the following result.
\begin{theorem} \label{thl32}
 Let $p\geq 2$. Assume for a constant $A,$
 $${\mathbf E} \Big[|\varepsilon_i|^p \, \Big| \, \sigma\{\varepsilon_j
, j\leq i-1\}  \Big] \leq A \ $$
for all $i \in [1, n].$ Then  for all  $x > 0,$
\begin{eqnarray}\label{indes10}
 {\mathbf P}\bigg( \pm (\theta_n -\theta)\sqrt{\Sigma_{k=1}^n \phi_{k}^2} \geq x   \bigg)   \   \leq \  \exp\Bigg\{- \frac{ x^2}{2 \, V^2 } \Bigg \} +  \frac{ C_p}{x^p}   \, ,\ \
\end{eqnarray}
where
\begin{eqnarray}\label{vinsa}
V^2=\frac{1}{4}(p+2)^2 e^p\, A^{2/p}    \ \ \   \textrm{and}\ \ \  C_p = \Big( 1+ \frac{2}{p} \Big)^p  A.
\end{eqnarray}
In particular, it holds  for any $x>0,$
\begin{eqnarray} \label{dfd1}
 {\mathbf P}\bigg( \pm (\theta_n -\theta)\sqrt{\Sigma_{k=1}^n \phi_{k}^2} \geq \sqrt{n}\, x   \bigg)   \   = \  O\bigg( \frac{c_x}{n^{p/2}} \bigg),\ \ \ \ \ \ \ \  n\rightarrow \infty,
\end{eqnarray}
where $c_x>0$ does not depend on $n.$
\end{theorem}

A similar inequality can be obtained by applying the Fuk inequality (\ref{s106}) to the martingale difference sequence  (cf.\ (\ref{sdgvf2}) for the definition of $(\xi _i,\mathcal{F}_i)_{i=1,...,n}$).
The Fuk inequality implies that  for all  $x > 0,$
\begin{eqnarray} \label{indes112}
 {\mathbf P}\bigg(  \pm(\theta_n -\theta)\sqrt{\Sigma_{k=1}^n \phi_{k}^2} \geq x   \bigg)   \   \leq \  \exp\Bigg\{- \frac{ x^2}{2 nV^2 } \Bigg \} +  \frac{ nC_p}{x^p},
\end{eqnarray}
where $V^2$ and $C_p$ are defined by (\ref{vinsa}).
In particular, it implies that for any $x>0,$
\begin{eqnarray} \label{dfd2}
 {\mathbf P}\bigg( \pm (\theta_n -\theta)\sqrt{\Sigma_{k=1}^n \phi_{k}^2} \geq \sqrt{n}\,  x   \bigg)   \   = \  O\bigg( \frac{c_x}{n^{p/2-1}} \bigg),\ \ \ \ \ n\rightarrow \infty,
\end{eqnarray}
where $c_x>0$ does not depend on $n.$
The order of  (\ref{dfd1}) is much better than that of (\ref{dfd2}).
Thus the refinement of (\ref{indes10}) on (\ref{indes112})  is significant.

If $(\varepsilon_k)_{k=1,...,n}$ have finite moments and uniformly bounded conditional variances,
by Theorem \ref{th23}, we obtain the following result which has the same polynomially decaying rate of Theorem  \ref{thl32}.
\begin{theorem} \label{thll3}
Let $p\geq 2.$ Assume for two constants $A$ and $B,$
$${\mathbf E} \Big[\varepsilon_i^2 \, \Big| \, \sigma\{\varepsilon_j
, j\leq i-1\} \Big]  \leq  A \  \ \ \ \ \ \textrm{and}\ \ \ \ \ \   {\mathbf E} \Big[|\varepsilon_i|^{p+\delta} \Big]  \leq B $$   for a small $\delta>0$ and all $i \in [1, n]$.
Then  for all $x  > 0$,
\begin{eqnarray}\label{fineq43}
 {\mathbf P}\bigg( \pm (\theta_n -\theta)\sqrt{\Sigma_{k=1}^n \phi_{k}^2} \geq x   \bigg)    \ \leq \ \exp\left\{- \frac{x^2}{2 \Big(A+ \frac{1}{3}x^{(2p+\delta)/(p+\delta) } \Big)} \right\}+  \frac{B}{  x^{p}} .
\end{eqnarray}
In particular, equality (\ref{dfd1}) holds.
\end{theorem}

In the following theorem, we assume that $(\varepsilon_i)_{i=1,...,n}$ have only a moment of order $p \in [1,2]$.
\begin{theorem}\label{thlecb}
Let $p\in [1, 2].$ Assume for a constant  $A$,
$${\mathbf E} \big[|\varepsilon_i|^p    \big]  \leq  A  $$   for  all $i \in [1, n]$.
Then  for all $x  > 0$,
\begin{eqnarray}
 {\mathbf P}\bigg( \pm (\theta_n -\theta)\sqrt{\Sigma_{k=1}^n \phi_{k}^2} \geq x   \bigg)    \ \leq \   \frac{2 A}{x^p} .
\end{eqnarray}
In particular, equality (\ref{dfd1}) holds.
\end{theorem}

%\begin{remark}
%The constant $2^{2-p}$ in (\ref{fineq43})  can be replaced by the more precise constant $\tilde C_p$ described in Proposition 1.8 of Pinelis \cite{P10}.
%\end{remark}

Theorems  \ref{thlin} and \ref{thlecb} focus on obtaining the large deviation inequalities. These inequalities do not depend on the distribution of
input random variables $(\phi_k)_{k=1,...,n}.$ Similar bounds are also expected to be  obtained  via the decoupling techniques of De la Pe\~{n}a \cite{D99} and  De la Pe\~{n}a and Gin\'{e} \cite{DG99}.  In particular, if $(\varepsilon_{k})_{k=1,...,n}$ are independent (instead of martingale differences),
with the method  of conditionally independent in De la Pe\~{n}a and Gin\'{e} \cite{DG99}, more precise bounds,  but depend  on the distribution of
input random variables, are allowed to be established.

Haeusler and Joos \cite{HJ88} proved that if the martingale differences satisfy  $\mathbf{E}[ |\xi_i|^{2+\delta}] < \infty$ for a constant $\delta > 0$ and all $i \in [1,n],$ then there exists a constant $C_\delta$, depending only on $\delta$, such that  for all $x \in \mathbf{R},$
\begin{eqnarray}
\bigg|\mathbf{P}(S_n\leq x) - \Phi(x) \bigg| \ \leq \ C_\delta \ \bigg( \sum_{i=1}^n\mathbf{E}\Big[  |\xi_i|^{2+\delta} \Big] +  \mathbf{E}\Big[ |\langle S\rangle_{n}-1|^{1+\delta/2} \Big] \bigg)^{1/(3+\delta)} \frac{1}{\displaystyle 1+|x|^{2+\delta}},\label{cdsi}
\end{eqnarray}
where $\Phi(x)=\frac{1}{\sqrt{2\pi}}\int_{-\infty}^{x}\exp\{-t^2/2\}dt$ is the standard normal distribution; see also  Hall and Heyde \cite{HH80} with the larger factor $\frac{1}{\displaystyle   1+|x|^{4(1+\delta/2)^2/(3+\delta)}}$ replacing $\frac{1}{\displaystyle 1+|x|^{2+\delta}}$.
Using (\ref{cdsi}), we obtain the following nonuniform Berry-Esseen bound, which depends on the distribution of
input random variables.

\begin{theorem}\label{thl33}
 Let $p> 2$. Assume that $(\varepsilon_i)_{i=1,...,n}$ satisfy ${\mathbf E} \big[ \varepsilon_i ^2 \ \big| \sigma\{\varepsilon_j
, j\leq i-1\}  \big] = \sigma^2 $ a.s.\ for a positive constant $\sigma$ and all $i \in [1, n].$
 Assume ${\mathbf E} [|\varepsilon_i|^p] \leq A$ for a constant $A$ and all $i \in [1, n]$.
Then   for all $x \in \mathbf{R},$
\begin{eqnarray}\label{sfsdvg}
 \bigg|{\mathbf P}\bigg(  (\theta_n -\theta)\sqrt{\Sigma_{k=1}^n \phi_{k}^2} \leq x \sigma  \bigg)   - \Phi(x) \bigg|  \  \leq \ C_p \, \Bigg( \sum_{i=1}^n  \mathbf{E}\Bigg[  \bigg| \frac{ \phi_{i}  }{ \sqrt{\Sigma_{k=1}^n \phi_{k}^2}} \bigg|^{p} \Bigg] \Bigg)^{1/(1+p)}  \frac{1}{1+|x|^{p}} \, ,
\end{eqnarray}
where $C_p$ is a constant depending only on $A, \sigma$ and $p$.
\end{theorem}

Notice that
\[
\sum_{i=1}^n  \mathbf{E}\Bigg[  \bigg| \frac{ \phi_{i}  }{ \sqrt{\Sigma_{k=1}^n \phi_{k}^2}} \bigg|^{p} \Bigg] \leq \sum_{i=1}^n  \mathbf{E}\Bigg[  \bigg| \frac{ \phi_{i}  }{ \sqrt{\Sigma_{k=1}^n \phi_{k}^2}} \bigg|^{2} \Bigg] =1.
\]
Thus (\ref{sfsdvg}) implies that the tail probability ${\mathbf P}\Big(  (\theta_n -\theta)\sqrt{\Sigma_{k=1}^n \phi_{k}^2} \geq x   \Big)$ has
the decaying rate $x^{-p}$ as $x\rightarrow \infty,$ which is coincident with the inequalities (\ref{indes10}) and (\ref{fineq43}).

\subsection{Weak invariance principles}\label{sec3.2}
In this subsection, let $(\xi _i,\mathcal{F}_i)_{i\geq 1}$ be  a
sequence of  stationary martingale differences. We have the following   weak
invariance principle for martingales.

The following  rate of convergence in the central limit theorem (CLT) for martingale difference sequences is due to
Ouchti (cf.\ Corollary 1 of \cite{O05}).  Assume that there exists a  constant $M >0$ such that
 $ \mathbf{E} [|\xi_i|^3 |\mathcal{F}_{i-1} ] \leq M  \mathbf{E} [\xi_i^2 |\mathcal{F}_{i-1} ]$
 a.s.\ for all $i \in \mathbf{N}.$ If the series $\sum_{i=1}^{\infty}  \mathbf{E} [\xi_i^2 |\mathcal{F}_{i-1} ]$  diverges a.s.\ and then there is a constant $C_M >0,$ depending on $M$, such that
\begin{eqnarray}
 \sup_{x \in \mathbf{R}}\bigg|{\mathbf P}\Big(  S_{v(n)} \leq x  \sqrt{n}  \Big)   - \Phi(x) \bigg|  \  \leq \ \frac{C_M}{ n^{1/4} } \,,
\end{eqnarray}
where
\[
v(n)  = \inf \Big\{ k \in \mathbf{N},  \ \ \  \left\langle S\right\rangle _k   \geq n \Big\} .
\]
Let
\[
H_n(t)= \frac{1}{\sqrt{ n  }} S_{v(\lfloor nt\rfloor)}  \ \ \ \ \ \  \ \textrm{for}\ 0\leq t \leq 1 .
\]
By Theorem \ref{th22}, we obtain
the following   weak
invariance principle for  martingales.
\begin{theorem} \label{f2jka}
Assume that there exists a  constant $M >0$ such that
 $ \mathbf{E} [|\xi_i|^3 |\mathcal{F}_{i-1} ] \leq M  \mathbf{E} [\xi_i^2 |\mathcal{F}_{i-1} ]$
 a.s.\ for all $i \in \mathbf{N}.$ If the series $\sum_{i=1}^{\infty}  \mathbf{E} [\xi_i^2 |\mathcal{F}_{i-1} ]$  diverges a.s.,
then the sequence of processes $\{H_n(t), 0\leq t \leq 1\}$ converges
in distribution to the standard Wiener process.
\end{theorem}

%%%%%%%%%%%%%%%%%%%%%%%%%%%%%%%%%%%%%%%%%%%%%%%%%%%%%%%%%%%%%%%%%%%%%%%%%%%%%%%%%%%%%%%%%%%%%%%%%%%%%%%%%%%%%%%%%%%%%%%%%%%%%%%%%%%%%%%%%%%%%%%%%%%%%%%%%%%%%

\section{Proof  of Theorem \ref{th21}}\label{sec4}
To prove Theorem \ref{th21},  we need the following technical lemma based on a truncation argument.
\begin{lemma}\label{lemma1} Assume   $\mathbf{E}[ \xi_i^2 \exp\{ |\xi_{i}|^{\alpha}  \}] < \infty$ for a constant $\alpha \in (0,1)$.
Set $\eta_{i}=\xi_{i}\mathbf{1}_{\{\xi_{i}\leq y\}}$ for  $y>0$.
 Then   for all $\lambda >0$,
\begin{eqnarray*}
\mathbf{E}[e^{\lambda \eta_i}|\mathcal{F}_{i-1}] \leq 1+  \frac{\lambda^2}{2} \mathbf{E}[ \eta_i^2 \exp\{\lambda  y^{1-\alpha}  ( \eta_{i} ^+)^{\alpha}  \}|\mathcal{F}_{i-1}].
\end{eqnarray*}
\end{lemma}

The proof of Lemma \ref{lemma1}  can be found in the proof of Proposition 3.5 in Dedecker  and Fan \cite{DF15}. However, instead of using the tower property of conditional expectation as in Dedecker  and Fan \cite{DF15}, we  use  changes of probability measure  in the proof of this theorem.
Set $\eta_{i}=\xi_{i}\mathbf{1}_{\{\xi_{i}\leq y\}}$ for some  $y>0.$ The exact value of $y$ is  given later.
Then  $(\eta_i,  \mathcal{F}_{i})_{i=1,...,n}$ is  a sequence of supermartingale differences, and it holds $ \mathbf{E}[\exp\left\{\lambda \eta_i   \right\} ] < \infty$ for all  $ \lambda \in (0, \infty)$ and all $i$. Define the exponential multiplicative
martingale  $Z(\lambda )=(Z_k(\lambda ),\mathcal{F}_k)_{k=0,...,n},$ where
\[
   Z_k(\lambda )=\prod_{i=1}^{  k}\frac{\exp\left\{\lambda \eta_i     \right\}}{\mathbf{E}\left[\exp\left\{\lambda \eta_i   \right\} | \mathcal{F}_{i-1} \right]},  \quad \quad  \quad Z_0(\lambda )=1. \label{C-1}
\]
If $T$ is a stopping time, then $Z_{T\wedge k}(\lambda )$ is also a martingale, where
\[
Z_{T\wedge k}(\lambda )=\prod_{i=1}^{T\wedge k}\frac{\exp\left\{\lambda \eta_i     \right\}}{\mathbf{E}\left[\exp\left\{\lambda \eta_i   \right\} |
\mathcal{F}_{i-1} \right]}, \quad \quad Z_0(\lambda )=1.
\]
Thus,  the random variable $Z_{T\wedge k}(\lambda ) $ is a
probability density on $(\Omega ,\mathcal{F},\mathbf{P})$, i.e.
$$ \int Z_{T\wedge k}(\lambda)  d \mathbf{P} = \mathbf{E}[Z_{T\wedge k}(\lambda)]=1.$$
Define the \emph{conjugate probability measure}
\begin{equation}
d\mathbf{P}_\lambda =Z_{T\wedge n}(\lambda )d\mathbf{P},   \label{chmeasure3}
\end{equation}
and denote by $\mathbf{E}_{\lambda}$ the expectation with
respect to $\mathbf{P}_{\lambda}.$
Since $\xi_i= \eta_i +\xi_{i}\mathbf{1}_{\{\xi_{i}> y\}},$ it follows that for any $x, y, u>0$,
\begin{eqnarray}
 &&  \mathbf{P}\Big(  S_k \geq x \ \mbox{and}\ \Upsilon(S)_k \leq u \ \mbox{for some}\ k \in [1,n] \Big)  \nonumber \\
 &\leq&  \mathbf{P}\left(   \sum_{i=1}^{k}\eta_i \geq x \ \mbox{and}\ \Upsilon(S)_k \leq u \ \mbox{for some}\ k \in [1,n] \right) \nonumber\\
  & & +\  \mathbf{P}\left(   \sum_{i=1}^{k}\xi_{i}\mathbf{1}_{\{\xi_{i}> y\}} > 0 \ \mbox{for some}\ k \in [1,n] \right)  \nonumber\\
&=:&  P_1 + \mathbf{P}\left( \max_{ 1\leq i\leq n}\xi_{i} > y \right). \label{fmuetoi}
\end{eqnarray}
For any $x, u>0$, define the stopping time
\[
T(x,u)=\min\bigg\{k\in [1, n]: \sum_{i=1}^{k}\eta_i \geq x \ \mbox{and}\ \Upsilon(S)_k \leq u \bigg\},
\]
with the convention that $\min{\emptyset}=0$. Then
 \[
\textbf{1}_{\{  S_k \geq x \ \mbox{and}\ \Upsilon(S)_k \leq u \ \mbox{for some}\ k \in [1,n]  \}} = \sum_{k=1}^{n}  \textbf{1}_{\{ T(x,u)=k\}}.
\]
By the change of measure (\ref{chmeasure3}), we deduce that for any $x,\lambda, u>0$,
\begin{eqnarray}
  P_1 &=& \mathbf{E}_{\lambda} \Big[ Z_{T\wedge n}(\lambda)^{-1}\textbf{1}_{\{S_k \geq x \ \mbox{and}\ \Upsilon(S)_k \leq u \ \mbox{for some}\ k \in [1,n]\}} \Big]\nonumber \\
 &=& \sum_{k=1}^{n}\mathbf{E}_{\lambda}\Big[ \exp\Big\{-\lambda \Big(\sum_{i=1}^{k}\eta_i \Big)+ \Psi_{k}(\lambda) \Big\} \textbf{1}_{\{T(x,u)=k\}}\Big], \label{ghnda}
\end{eqnarray}
where
\begin{eqnarray}\label{fsbdphi}
\Psi_{k}(\lambda)= \sum_{i=1}^k \log \mathbf{E} \left[\exp\left\{\lambda \eta_i   \right\} |
\mathcal{F}_{i-1} \right].
\end{eqnarray}
Set $\lambda=y^{\alpha-1}.$ By Lemma \ref{lemma1} and the inequality $\log(1+t) \leq t$ for all $t\geq 0$,  it is easy to see that for any  $x>0$,
\begin{eqnarray*}
\Psi_{k}(\lambda) &\leq& \sum_{i=1}^k  \log \bigg( 1+  \frac{\lambda^2}{2} \mathbf{E}[ \eta_i^2 \exp\{\lambda  y^{1-\alpha}  ( \eta_{i} ^+)^{\alpha}  \}|\mathcal{F}_{i-1}] \bigg)  \\
&\leq& \sum_{i=1}^k  \frac{\lambda^2}{2} \mathbf{E}[ \eta_i^2 \exp\{\lambda  y^{1-\alpha}  ( \eta_{i} ^+)^{\alpha}  \}|\mathcal{F}_{i-1}]\\
&\leq& \frac12  y^{2\alpha-2} \Upsilon(S)_k.
\end{eqnarray*}
By the fact that $\sum_{i=1}^{k}\eta_i \geq x$  and  $\Psi_{k}(\lambda) \leq  \frac12  y^{2\alpha-2} u$ on the set  $\{T(x,u)=k\}$. we find that
for any $x,  u>0$,
\begin{eqnarray*}
 P_1   &\leq&  \exp \left\{ - \lambda x +   \frac12  y^{2\alpha-2} u \right\} \mathbf{E}_{\lambda} \Big[ \sum_{k=1}^{n}\textbf{1}_{\{T(x,u)=k\}} \Big] \nonumber\\
 &\leq&  \exp \left\{ -  y^{\alpha-1} x +   \frac12  y^{2\alpha-2} u   \right\} .
\end{eqnarray*}
From (\ref{fmuetoi}), it follows that
\begin{eqnarray}
 && \mathbf{P}\Big(  S_k \geq x \ \mbox{and}\ \Upsilon(S)_k \leq u \ \mbox{for some}\ k \in [1,n] \Big)  \nonumber \\
 &&\ \ \ \ \ \ \ \ \ \ \ \ \ \ \ \ \ \ \ \leq\
 \exp \left\{ -  y^{\alpha-1} x +   \frac12  y^{2\alpha-2} u  \right\} + \mathbf{P}\left( \max_{ 1\leq i\leq n}\xi_{i} > y \right).  \label{sdfo}
\end{eqnarray}
By the exponential Markov  inequality, we have the following estimation: for any $x>0$,
\begin{eqnarray}
\mathbf{P}\left( \max_{ 1\leq i\leq n}\xi_{i} > y \right) &\leq& \sum_{i=1}^n \mathbf{P}\left(  \xi_{i} > y \right) \nonumber\\
 &\leq& \frac{1}{y^{2}} \exp\{ -y^{\alpha} \} \ \sum_{i=1}^n \mathbf{E} [ \xi_i^2 \exp\{ ( \xi_{i} ^+)^{\alpha} \} ] \nonumber\\
 &\leq&\frac{ C_n}{y^2}\exp \left\{- y^{\alpha}   \right\}. \label{ft90}
\end{eqnarray}
Taking
\begin{eqnarray*}
y = \left\{ \begin{array}{ll}
\Big(\displaystyle \frac u x\Big)^{1/(1-\alpha)} & \textrm{\ \ \ \ \ if $0\leq x < u^{1/(2-\alpha)}$}  \\
\vspace{-0.3cm}\\
x  & \textrm{\ \ \ \ \ if $x \geq u^{1/(2-\alpha)}$, }
\end{array} \right.
\end{eqnarray*}
from  (\ref{sdfo}) and (\ref{ft90}),    we obtain the desired inequality.
This completes the proof of  Theorem \ref{th21}.
\hfill\qed

\vspace{0.3cm}

\noindent\emph{Proof of Corollary \ref{th1}.} Set $u_n=||\Upsilon(S)_n||_{\infty}.$ Then $u_n =o(n^{2-\alpha}),  n \rightarrow \infty,$ by the assumptions of Theorem \ref{th1}. For any $x>0,$ by Theorem \ref{th21}, we have
\begin{eqnarray*}
  \mathbf{P}\left( \max_{1\leq k \leq n}  S_k   \geq nx \right)    &\leq& \exp \left\{ - (nx)^\alpha \left( 1- \frac{u_n}{2\, (nx)^{2-\alpha}}\right)  \right\} + \frac{ C_n}{(nx)^2}\exp \Big\{- (nx)^{\alpha}   \Big\} \\
  \\
  &\leq& \Big( 1+ \frac{ C_n}{(nx)^2} \Big) \exp \left\{ - (nx)^\alpha \left( 1- \frac{u_n}{2\, (nx)^{2-\alpha}}\right)  \right\} .
\end{eqnarray*}
Since $u_n  \geq C_n$, we have  $C_n =o(n^{2-\alpha}),  n \rightarrow \infty.$
Hence it holds
\begin{eqnarray*}
 \limsup_{n\rightarrow \infty}\frac{1}{n^\alpha}\log \mathbf{P}\left( \max_{1\leq k \leq n}  S_k     \geq nx \right) \ \leq  \ - x ^\alpha \, .
\end{eqnarray*}
This completes the proof of  Corollary \ref{th1}.
\hfill\qed

\vspace{0.3cm}

\noindent\emph{Proof of Remark \ref{re1}.}
 Note that
\begin{eqnarray*}
 \mathbf{E}[ \xi_1^2 \exp\{ |\xi_{1}|^{\alpha} \}]
&=& \int_0^\infty \mathbf{P}(|\xi_1|\geq x) \Big( 2x + \alpha x^{1+\alpha} \Big)e^{x^\alpha} dx \  <\  \infty.
\end{eqnarray*}
Thus
\begin{eqnarray*}
\Upsilon(S)_n \ \leq \ n \mathbf{E}[ \xi_1^2 \exp\{ |\xi_{1}|^{\alpha} \}] \ =\ o(n^{2-\alpha}),\ \ \ \ \ \ \ n \rightarrow \infty.
\end{eqnarray*}
It is easy to see that  for any $x, \varepsilon> 0,$ we have
\begin{eqnarray*}
 \mathbf{P}\left( \max_{1\leq k \leq n} S_k \geq nx \right) &\geq&  \mathbf{P}\Big(   S_n \geq nx \Big) \\
 &\geq&  \mathbf{P}\bigg(\sum_{i=2}^n \xi_i \geq -n\varepsilon, \xi_1 \geq n(x+\varepsilon)\bigg) \\
 &=&\mathbf{P} \Big( \sum_{i=2}^n \xi_i \geq -n\varepsilon \Big) \ \mathbf{P}\Big(  \xi_1 \geq n(x+\varepsilon)\Big).
\end{eqnarray*}
The first probability on the right-hand side trends to $1$ as $n\rightarrow \infty$ due to the law of large numbers.
By (\ref{as11}), the second term on the right-hand side has the following lower bound
\begin{eqnarray*}
\mathbf{P} \Big(  \xi_1 \geq n(x+\varepsilon) \Big)  \ \geq \ \Big(n(x+\varepsilon)\Big)^{-2p} \, \exp \bigg\{ - \Big(n(x+\varepsilon)\Big)^\alpha  \bigg\}
\end{eqnarray*}
for all $n$ large enough.  Hence
\begin{eqnarray*}
 \lim_{n\rightarrow \infty}\frac{1}{n^\alpha}\log \mathbf{P}\left( \max_{1\leq k \leq n} \frac{1}{n}S_k \geq x \right) \ \geq\ -  (x+\varepsilon) ^\alpha . \,
\end{eqnarray*}
Letting $\varepsilon \rightarrow 0,$  we obtain
\begin{eqnarray*}
 \lim_{n\rightarrow \infty}\frac{1}{n^\alpha}\log \mathbf{P}\left( \max_{1\leq k \leq n} \frac{1}{n}S_k \geq x \right) \ \geq\ -   x^\alpha . \,
\end{eqnarray*}
Combining this result with Theorem \ref{th1}, we get (\ref{fineq15}).\hfill\qed

\section{Proof of Theorem \ref{th22}}\label{sec5}
To prove Theorem \ref{th22},  we need the following technical lemma.
\begin{lemma}\label{lem1}
 Let $p\geq 2$. Assume ${\mathbf E} [|\xi_i|^p]   < \infty $ for all $i \in [1, n]$.
Set $\eta_{i}=\xi_{i}\mathbf{1}_{\{\xi_{i}\leq y\}}$ for  $y>0$.
 Then   for all $\lambda >0$,
\begin{eqnarray*}
\mathbf{E}[e^{\lambda \eta_i}|\mathcal{F}_{i-1}] \leq   1+  \frac{1}{2}e^p\lambda^2 \mathbf{E}[ \xi_i^2  |\mathcal{F}_{i-1}] +f(y)\mathbf{E} [(\xi_i^+)^p |\mathcal{F}_{i-1} ] \, ,
\end{eqnarray*}
where the function
\begin{eqnarray}\label{funcfx}
f(u)=\frac{e^{\lambda u }-1-\lambda u}{u^p}, \ \ \ u> 0.
\end{eqnarray}
\end{lemma}

\vspace{0.3cm}

\noindent\emph{Proof.} We argue as in Fuk and Nagaev \cite{FN71} (see also Fuk \cite{Fu73}). Using a two term Taylor's expansion, we have for some $\theta \in [0, 1],$
\begin{eqnarray*}
 e^{\lambda \eta_i}  &\leq& 1+ \lambda \eta_i +  \frac{\lambda^2}{2} \eta_i^2 \mathbf{1}_{\{\lambda \eta_i \leq p  \}}  e^{  \lambda \theta \eta_i  } + f(\eta_i) (\eta_i^+)^p  \mathbf{1}_{\{\lambda \eta_i > p  \}} .
\end{eqnarray*}
Remark that the function $f$ is positive and increasing for $\lambda u \geq p.$
Since $\mathbf{E}[ \eta_i |\mathcal{F}_{i-1}]\leq \mathbf{E}[\xi_i |\mathcal{F}_{i-1}] = 0$ and $\eta_{i} \leq y$, it follows that
\begin{eqnarray*}
\mathbf{E}[e^{\lambda \eta_i}|\mathcal{F}_{i-1}]  &\leq & 1+  \frac{1}{2}e^p\lambda^2 \mathbf{E}[ \eta_i^2  |\mathcal{F}_{i-1}] +f(y)\mathbf{E} [(\eta_i^+)^p |\mathcal{F}_{i-1} ] \\
 &\leq & 1+  \frac{1}{2}e^p\lambda^2 \mathbf{E}[ \xi_i^2  |\mathcal{F}_{i-1}] +f(y)\mathbf{E} [(\xi_i^+)^p |\mathcal{F}_{i-1} ] \, ,
\end{eqnarray*}
which gives the desired inequality.\hfill\qed
\vspace{0.3cm}

 We  make use of Lemma \ref{lem1} to prove Theorem \ref{th22}.
Set $\eta_{i}=\xi_{i}\mathbf{1}_{\{\xi_{i}\leq y\}}$ for   $y>0$.
Define the  conjugate probability measure  $d\mathbf{P}_\lambda$ by (\ref{chmeasure3})
and denote by $\mathbf{E}_{\lambda}$ the expectation with
respect to $\mathbf{P}_{\lambda}.$
Since $\xi_i= \eta_i +\xi_{i}\mathbf{1}_{\{\xi_{i}> y\}},$ it follows that for any $x, y, u, w>0$,
\begin{eqnarray}
 &&  \mathbf{P}\left(   S_k >x,\ \langle S\rangle_k \leq v  \  \textrm{and} \ \Xi(S)_k   \leq  w \textrm{ for some }k \in [1, n]\right)  \nonumber \\
 &\leq&  \mathbf{P}\left(   \sum_{i=1}^{k}\eta_i \geq x,  \ \langle S\rangle_k \leq v  \  \textrm{and} \ \Xi(S)_k   \leq  w \textrm{ for some }k \in [1, n] \right) \nonumber\\
  & & +\  \mathbf{P}\left(   \sum_{i=1}^{k}\xi_{i}\mathbf{1}_{\{\xi_{i}> y\}} > 0 \ \mbox{for some}\ k \in [1,n] \right)  \nonumber\\
&=:&  P_2 + \mathbf{P}\left( \max_{ 1\leq i\leq n}\xi_{i} > y \right).   \label{fmstdfoi}
\end{eqnarray}
For any $x, v, w>0$, define the stopping time $T:$
\[
T(x,v,w)=\min\bigg\{k\in [1, n]:  S_k \geq x,\ \langle S\rangle_k \leq v  \  \textrm{and} \ \Xi(S)_k   \leq  w \bigg\},
\]
with the convention that $\min{\emptyset}=0$. Then
 \[
\textbf{1}_{\{  S_k >x,\ \langle S\rangle_k \leq v  \  \textrm{and} \ \Xi(S)_k   \leq  w \textrm{ for some }k \in [1, n]  \}} = \sum_{k=1}^{n}  \textbf{1}_{\{ T =k\}}.
\]
By the change of measure (\ref{chmeasure3}), we deduce that for any $x, y, \lambda, u, w>0$,
\begin{eqnarray*}
  P_2 &=& \mathbf{E}_{\lambda} \Big[ Z_{T\wedge n}(\lambda)^{-1}\textbf{1}_{\{ S_k >x,\ \langle S\rangle_k \leq v  \  \textrm{and} \ \Xi(S)_k   \leq  w \textrm{ for some }k \in [1, n] \}} \Big]\nonumber \\
 &=& \sum_{k=1}^{n}\mathbf{E}_{\lambda}\Big[ \exp\Big\{-\lambda \Big(\sum_{i=1}^{k}\eta_i \Big)+ \Psi_{k}(\lambda) \Big\} \textbf{1}_{\{T=k\}}\Big],
\end{eqnarray*}
where $\Psi_{k}(\lambda)$ is defined by (\ref{fsbdphi}).
By Lemma \ref{lem1} and the inequality $\log(1+t) \leq t$ for $t\geq 0$,  it is easy to see that for any  $x, y, \lambda, u, w>0$,
\begin{eqnarray*}
\Psi_{k}(\lambda) &\leq& \sum_{i=1}^k  \log \bigg(  1+  \frac{1}{2}e^p\lambda^2 \mathbf{E}[ \xi_i^2  |\mathcal{F}_{i-1}] +f(y)\mathbf{E} [(\xi_i^+)^p |\mathcal{F}_{i-1} ] \bigg)  \\
&\leq& \sum_{i=1}^k  \bigg(  \frac{1}{2}e^p\lambda^2 \mathbf{E}[ \xi_i^2  |\mathcal{F}_{i-1}] +f(y)\mathbf{E} [(\xi_i^+)^p |\mathcal{F}_{i-1} ] \bigg),
\end{eqnarray*}
where $f(y)$ is defined by (\ref{funcfx}).
By the fact that $\sum_{i=1}^{k}\eta_i \geq x$  and  $\Psi_{k}(\lambda) \leq   \frac{1}{2}e^p\lambda^2 v +f(y)w$ on the set  $\{T=k\}$. we find that
for any $x, y, \lambda, u, w>0$,
\begin{eqnarray*}
 P_2   &\leq&  \exp \left\{ - \lambda x +  \frac{1}{2}e^p\lambda^2 v +f(y)w \right\} \mathbf{E}_{\lambda} \Big[ \sum_{k=1}^{n}\textbf{1}_{\{T=k\}} \Big] \nonumber\\
 &\leq&  \exp \left\{  - \lambda x + \frac{1}{2}e^p\lambda^2 v +f(y)w  \right\} .
\end{eqnarray*}
Next we carry out an argument as in  Fuk and Nagaev \cite{FN71}.
Then
\begin{eqnarray}\label{endineq}
 P_2   &\leq&   \exp\Bigg\{- \frac{\alpha^2x^2}{2e^p v} \Bigg \} + \exp\Bigg\{- \frac{\beta x}{y} \log\Big( 1 + \frac{\beta x y^{p-1}}{w}\Big)  \Bigg \} \, ,
\end{eqnarray}
where $\alpha$ and $\beta$ are defined by (\ref{defab}).
Combining the inequalities (\ref{fmstdfoi}) and (\ref{endineq}) together, we obtain the desired inequality.
This completes the proof of Theorem \ref{th22}
\hfill\qed

\vspace{0.3cm}

\noindent\emph{Proof of Corollary \ref{co22}.} When $y=\beta x,$ from (\ref{sdfds}), it is easy to see that for all $x> 0,$
\begin{eqnarray*}
\mathbf{P}\left( \max_{ 1\leq i\leq n}\xi_{i} > y \right) &\leq& \sum_{i=1}^n\mathbf{P}\Big(  \xi_{i} > \beta x \Big)
\ \leq\ \frac{1}{\beta^p x^p }  \sum_{i=1}^n\mathbf{E}[|\xi_{i}|^p]
\ \leq \  \frac{C_p}{ x^p }
\end{eqnarray*}
and
\begin{eqnarray*}
 \exp\Bigg\{- \frac{\beta x}{y} \log\Big( 1 + \frac{\beta x y^{p-1}}{w}\Big)  \Bigg \}\ \leq\ \frac{w}{\beta x y^{p-1} + w}
\ \leq\ \frac{w}{\beta^p x^p }\ \leq \  \frac{C_p}{ x^p } ,
\end{eqnarray*}
where $C_p$ is defined by (\ref{cpdef}). Thus (\ref{sdfds}) implies (\ref{sfddghfd}).\hfill\qed

\section{Proofs of Theorem \ref{th23} and Corollary \ref{co23}}
To prove Theorem \ref{th23}, we need the following inequality whose proof can be found in Fan, Grama and Liu \cite{FGL12} (cf. Corollary 2.3 and Remark 2.1 therein).
\begin{lemma}
\label{leedm1} Assume ${\mathbf E} [\xi_i^2]   < \infty $ for all $i \in [1, n]$.  Then
for all $ x, y, v > 0$,
\begin{eqnarray}
&&{\mathbf P}\Big( S_k \geq x  \  \textrm{and} \ \langle S\rangle_k \leq v^2   \textrm{ for some }k \in [1, n] \Big) \ \ \ \ \ \ \ \ \ \nonumber \\
 &&\ \ \ \  \ \ \ \ \ \ \ \ \ \ \ \leq \  \exp\left\{- \frac{x^2}{2(v^2+ \frac{1}{3}xy)} \right\}+ \mathbf{P}\left( \max_{1\leq i \leq n} \xi_{i}>y  \right). \nonumber
\end{eqnarray}
\end{lemma}
\noindent\emph{Proof of Theorem \ref{th23}.}  By Lemma \ref{leedm1}  and the Markov inequality,  it follows that for all $ x, y, v > 0$,
\begin{eqnarray}
&&{\mathbf P}\Big( S_k \geq x  \  \textrm{and} \ \langle S\rangle_k \leq v^2   \textrm{ for some }k \in [1, n] \Big) \ \ \ \ \ \ \ \ \ \nonumber \\
 &&\ \ \ \  \ \ \ \leq \  \exp\left\{- \frac{x^2}{2(v^2+ \frac{1}{3}xy)} \right\}+ \sum_{i=1}^n\mathbf{P}\Big(  \xi_{i}>y  \Big)  \nonumber \\
 &&\ \ \ \ \ \ \ \ \ \  \ \ \leq \ \exp\left\{- \frac{x^2}{2(v^2+ \frac{1}{3}xy)} \right\}+ \frac{1}{y^{p+\delta}} \sum_{i=1}^n{\mathbf E} \Big[|\xi_i|^{p+\delta}\mathbf{1}_{\{\xi_{i}>y\}}\Big]. \nonumber
\end{eqnarray}
Taking $y=x^{p/(p+\delta) }$ in the last inequality, we obtain the desired inequality. This completes the proof of Theorem \ref{th23}. \hfill\qed

\vspace{0.3cm}

\noindent\emph{Proof of Corollary \ref{co23}.} Notice that $p+\delta > 2.$ It is easy to see that  for any $x, v> 0,$
\begin{eqnarray*}
 {\mathbf P}\Big( \max_{1\leq k \leq n} S_k \geq x   \Big)  &\leq&  {\mathbf P}\Big( \max_{1\leq k \leq n} S_k \geq x \ \ \textrm{and} \ \langle S\rangle_n \leq nv^2   \Big) + {\mathbf P}\Big(  \langle S\rangle_n > nv^2   \Big) \\
  &\leq& {\mathbf P}\Big( S_k \geq x  \  \textrm{and} \ \langle S\rangle_k \leq n v^2   \textrm{ for some }k \in [1, n] \Big)+ {\mathbf P}\Big(  \langle S\rangle_n > nv^2   \Big)\\
  &\leq& {\mathbf P}\Big( S_k \geq x  \  \textrm{and} \ \langle S\rangle_k \leq n v^2   \textrm{ for some }k \in [1, n] \Big)\\
  && + \ \frac{\mathbf{E}[|\langle S\rangle_n|^{(p+\delta)/2}]}{n^{(p+\delta) /2} v^{p+\delta }} ,
\end{eqnarray*}
which gives the first desired inequality.
By the H\"{o}lder inequality,
it follows that
\[
\sum_{i=1}^{n} a_i \ \leq \ n^{1- 2/(p+\delta) } \Big(\sum_{i=1}^{n} a_i^{(p+\delta)/2} \Big)^{2/(p+\delta)},   \ \ \ \ \ \ a_i\geq 0, \ i=1,...,n.
\]
Hence
\[
\Big(\sum_{i=1}^{n} a_i \Big)^{(p+\delta)/2} \ \leq \  n^{(p-2+\delta)/2 }  \sum_{i=1}^{n} a_i^{(p+\delta)/2}  ,  \ \ \ \ \ \ \ a_i\geq 0, \ i=1,...,n.
\]
Then we have
\begin{eqnarray*}
 \mathbf{E}[|\langle S\rangle_n|^{(p+\delta)/2}]
   &\leq&  n^{(p-2+\delta)/2 } \sum_{i=1}^n {\mathbf E}\Big[  {\mathbf E} \big[ \xi_i^2|  \mathcal{F}_{i-1} \big]^{(p+\delta)/2}  \Big] \\
   &\leq&  n^{(p-2+\delta)/2 } \sum_{i=1}^n {\mathbf E}\Big[  {\mathbf E} \big[ |\xi_i|^{p+\delta} |  \mathcal{F}_{i-1} \big] \Big] \\
  & =&  n^{(p-2+\delta)/2 } \sum_{i=1}^n{\mathbf E} \big[|\xi_i|^{p+\delta} \big].
\end{eqnarray*}
This completes the proof of  corollary.\hfill\qed

%%%%%%%%%%%%%%%%%%%%%%%%%%%%%%%%%%%%%%%%%%%%%%%%%%%%%%%%%%%%%%%%%%%%%%%%%%%%%%%%%%%%%%%%%%%%%%%%%%%%%%%%%%%%%%%%%%%%%%%%%%%%%%%%%%%%%%%%%%%%%%%%%%%%%%%%%%%%%
\section{Proofs  of Theorems \ref{thlin} - \ref{thl33} } \label{sec6}
From (\ref{ine29}) and (\ref{ine30}), it is easy to see that
\begin{equation}\label{sdgvf1}
\theta_n -\theta = \sum_{k=1}^n\frac{ \phi_{k} \varepsilon_k}{\Sigma_{k=1}^n \phi_{k}^2}. \nonumber
\end{equation}
For any $i=1,...,n$, set
\begin{eqnarray}\label{sdgvf2}
\xi_i= \frac{ \phi_{i} \varepsilon_i}{ \sqrt{\Sigma_{k=1}^n \phi_{k}^2}}\ \ \ \  \textrm{and} \ \ \ \ \mathcal{F}_{i} = \sigma \Big\{ \phi_{k}, \varepsilon_k, 1\leq k\leq i,\  \phi_{k}^2, 1\leq k\leq n \Big\}.
\end{eqnarray}
Then $(\xi _i,\mathcal{F}_i)_{i=1,...,n}$ is a sequence of martingale differences, and satisfies
\begin{eqnarray}
 S_n= \sum_{i=1}^n\xi_i = (\theta_n -\theta)\sqrt{\Sigma_{k=1}^n \phi_{k}^2} \ .
\end{eqnarray}

\vspace{0.3cm}

\noindent\emph{Proof of Theorem \ref{thlin}.}
Notice that
$$  \Upsilon  (S)_n  \ \leq  \ \sum_{i=1}^{n} \frac{ \phi_{i}^2 }{   \Sigma_{k=1}^n \phi_{k}^2 } \mathbf{E}[\varepsilon_i^2\exp\{|\varepsilon_{i}|^{\alpha} \} | \mathcal{F}_{i-1} ] \ \leq \  \sum_{i=1}^{n} \frac{ \phi_{i}^2 D}{  \Sigma_{k=1}^n \phi_{k}^2 } \ = \ D . $$
Applying Theorem \ref{th21}  to $(\xi _i,\mathcal{F}_i)_{i=1,...,n},$  we find that (\ref{ineq6}), with $u\geq \max\{D, 1 \}$, is an upper bound on the tail probabilities $\mathbf{P}\left(  (\theta_n -\theta)\sqrt{ \Sigma _{k=1}^n \phi_{k}^2} \geq x \right).$

Similarly, applying Theorem \ref{th21}  to $(-\xi _i,\mathcal{F}_i)_{i=1,...,n},$ we find that (\ref{ineq6}), with $u\geq \max\{D, 1 \}$, is also an upper bound on the tail probabilities   $\mathbf{P}\left(  -(\theta_n -\theta)\sqrt{ \Sigma _{k=1}^n \phi_{k}^2} \geq x \right).$ This completes the proof   of  Theorem \ref{thlin}.\hfill\qed

\vspace{0.3cm}

\noindent\emph{Proof of Theorem \ref{thFA12}.}
By the fact $$\mathbf{E}[\varepsilon_i^2   | \mathcal{F}_{i-1} ] = \mathbf{E}[\varepsilon_i^2   | \sigma\{\varepsilon_k, 1\leq k \leq i-1  \}  ] \leq E ,$$ it follows that
$$  \langle S\rangle_n \ \leq \  \sum_{i=1}^{n} \frac{ \phi_{i}^2 }{  (\Sigma_{k=1}^n \phi_{k}^2)} \mathbf{E}[\varepsilon_i^2  | \mathcal{F}_{i-1} ] \ \leq \ E.  $$
Similarly, by the fact $ \mathbf{E}\big[ \exp\{|\varepsilon_i|^{\frac{\alpha}{1-\alpha}} \}   \big] \leq F,$
it is easy to see that
$$  \mathbf{E} [\exp\{|\xi_i|^{\frac{\alpha}{1-\alpha}} \}  ] \ \leq \ \mathbf{E}\big[ \exp\{|\varepsilon_i|^{\frac{\alpha}{1-\alpha}} \}   \big] \ \leq \ F.$$
Applying Theorem 2.2 of Fan, Grama and Liu \cite{Fx1} to $(\pm\xi _i,\mathcal{F}_i)_{i=1,...,n},$     we obtain  the desired inequality.\hfill\qed

\vspace{0.3cm}

\noindent\emph{Proof of Theorem \ref{thl32}.}
By the fact $$\mathbf{E}[|\varepsilon_i|^p   | \mathcal{F}_{i-1} ] = \mathbf{E}[|\varepsilon_i|^p   | \sigma\{\varepsilon_k, 1\leq k \leq i-1  \}  ] \leq A ,$$ it follows that
$$  \langle S\rangle_n \ \leq \  \sum_{i=1}^{n} \frac{ \phi_{i}^2 }{  (\Sigma_{k=1}^n \phi_{k}^2)} \mathbf{E}[\varepsilon_i^2  | \mathcal{F}_{i-1} ] \ \leq \ \sum_{i=1}^{n} \frac{ \phi_{i}^2  }{ (\Sigma_{k=1}^n \phi_{k}^2) } \Big(\mathbf{E}[|\varepsilon_i|^p   | \mathcal{F}_{i-1} ] \Big)^{2/p} = A^{2/p}   $$
and
$$   \sum_{i=1}^{n} \mathbf{E} [|\xi_i|^p |\mathcal{F}_{i-1} ] \ \leq \  \sum_{i=1}^{n} \frac{ \phi_{i}^2 }{  (\Sigma_{k=1}^n \phi_{k}^2)} \mathbf{E}[|\varepsilon_i|^p   | \mathcal{F}_{i-1} ] \ \leq \ A.$$
Applying Corollary \ref{co22} to $(\pm\xi _i,\mathcal{F}_i)_{i=1,...,n},$  we obtain the desired inequality.\hfill\qed

\vspace{0.3cm}

\noindent\emph{Proof of Theorem \ref{thll3}.}
By the fact $$\mathbf{E}[\varepsilon_i^2   | \mathcal{F}_{i-1} ] = \mathbf{E}[\varepsilon_i^2   | \sigma\{\varepsilon_k, 1\leq k \leq i-1  \}  ] \leq A ,$$
it follows that
$$  \langle S\rangle_n  \ \leq \  \sum_{i=1}^{n} \frac{ \phi_{i}^2 }{  (\Sigma_{k=1}^n \phi_{k}^2)} \mathbf{E}[\varepsilon_i^2  | \mathcal{F}_{i-1} ] \ \leq \ A .  $$
Similarly,  by the fact $\mathbf{E}[|\varepsilon_i|^{p+\delta}  ] \leq B ,$
it follows that
$$   \sum_{i=1}^{n} \mathbf{E} [|\xi_i|^{p+\delta} ] \ = \ \sum_{i=1}^{n} \mathbf{E} \bigg[ \Big| \frac{ \phi_{i}^2 }{   \Sigma_{k=1}^n \phi_{k}^2 } \Big|^{\frac{p+\delta}{2}} \bigg]  \mathbf{E} \big[|\varepsilon_i|^{p+\delta}\big] \ \leq \ \mathbf{E} \bigg[ \sum_{i=1}^{n}   \frac{ \phi_{i}^2 }{   \Sigma_{k=1}^n \phi_{k}^2 }     \bigg] B \ =\  B.$$
Applying Theorem \ref{th23} to $(\pm\xi _i,\mathcal{F}_i)_{i=1,...,n},$  we obtain the desired inequality. \hfill\qed

\vspace{0.3cm}

\noindent\emph{Proof of Theorem \ref{thlecb}.} Let $p \in [1, 2].$  By the inequality $$\Big(\sum_{i=1}^n a_i \Big)^\alpha \leq \sum_{i=1}^n a_i^\alpha,\ \ \ \  a_i \geq 0 \ \ \textrm{and}\ \  \alpha \in (0, 1],$$
we have
$$   \sum_{i=1}^{n} \mathbf{E} [|\xi_i|^{p } ] \ = \ \sum_{i=1}^{n} \mathbf{E} \bigg[   \frac{ (\phi_{i}^2 )^{p/2 }}{  ( \Sigma_{k=1}^n \phi_{k}^2 )^{p/2 } }   \bigg]  \mathbf{E} \big[|\varepsilon_i|^{p }\big]  \ \leq \  \mathbf{E} \bigg[   \frac{ (\sum_{i=1}^{n}\phi_{i}^2 )^{p/2 }}{  ( \Sigma_{k=1}^n \phi_{k}^2 )^{p/2 } }   \bigg]  A \ = \ A.$$
By the inequality of von Bahr and Esseen  (cf. Theorem 2 of \cite{V65}), we get
$$\mathbf{E} [ |S_n|^{p } ] \ \leq \ 2 \sum_{i=1}^{n} \mathbf{E} [|\xi_i|^{p } ] \ \leq\ 2A.$$
Then for all $x>0,$
\begin{eqnarray*}
 {\mathbf P}\bigg( \pm (\theta_n -\theta)\sqrt{\Sigma_{k=1}^n \phi_{k}^2} \geq x   \bigg) \ =\ {\mathbf P}\Big( \pm S_n \geq x   \Big)   \ \leq \   \frac{\mathbf{E} [ |S_n|^{p } ] }{x^p}  \ \leq \   \frac{2 A}{x^p}  .
\end{eqnarray*}
This completes the proof of theorem.\hfill\qed

\vspace{0.3cm}

\noindent\emph{Proof of Theorem \ref{thl33}.} It is obvious that
$$   \frac{\theta_n -\theta}{\sigma}  \sqrt{\Sigma_{k=1}^n \phi_{k}^2}  =\sum_{i=1}^n \eta_i,$$
where $\eta_i=\xi_i/\sigma.$
Notice that  $\mathbf{E}[\varepsilon_i^2  | \mathcal{F}_{i-1} ]={\mathbf E} \big[ \varepsilon_i ^2 \ \big| \sigma\{\varepsilon_j
, j\leq i-1\}  \big] = \sigma^2$ a.s..\ Then we have
$$  \sum_{i=1}^{n} \mathbf{E} [\eta_i^2 |\mathcal{F}_{i-1} ] \ = \ \frac{ \langle S\rangle_n}{\sigma^2} \ =\  \sum_{i=1}^{n} \frac{ \phi_{i}^2 }{  (\Sigma_{k=1}^n \phi_{k}^2)}\frac{ \mathbf{E}[\varepsilon_i^2  | \mathcal{F}_{i-1} ]}{\sigma^2} \ =\ \sum_{i=1}^{n} \frac{ \phi_{i}^2  }{  \Sigma_{k=1}^n \phi_{k}^2 }  \ = \ 1   $$
and
$$  \sum_{i=1}^{n} \mathbf{E} [|\eta_i|^p |\mathcal{F}_{i-1} ]  \ \leq \  \sum_{i=1}^{n}   \mathbf{E}\bigg[  \bigg| \frac{ \phi_{i}  }{ \sqrt{\Sigma_{k=1}^n \phi_{k}^2}} \bigg|^{p} \bigg] \frac{ \mathbf{E}[|\varepsilon_i|^p   | \mathcal{F}_{i-1} ]}{\sigma^p} \ \leq\ \frac{A}{\sigma^p}\sum_{i=1}^{n}   \mathbf{E}\bigg[  \bigg| \frac{ \phi_{i}  }{ \sqrt{\Sigma_{k=1}^n \phi_{k}^2}} \bigg|^{p} \bigg].$$
Applying inequality  (\ref{cdsi}) to the martingale difference sequence $(\eta_i,\mathcal{F}_i)_{i=1,...,n}$ with $\delta=p-2$,  we obtain the desired inequality.\hfill\qed

\section{Proof   of Theorem \ref{f2jka} } \label{sec7}
The proof is based on the result  of   Ouchti \cite{O05}.

\vspace{0.2cm}

\noindent\emph{Proof of Theorem \ref{f2jka}.} By the rate of convergence in the CLT for martingale difference sequences of
Ouchti (cf.\ Corollary 1 of \cite{O05}), it suffices to verify the tightness of $H_n.$  By Theorem 8.4 of Billingsley \cite{B68} for  stationary martingale difference sequences, we only need to show that for any $\varepsilon > 0,$ there exist a $\lambda,$ with $\lambda >1,$ and an integer $n_0$ such that for every
$n \geq n_0,$
\begin{equation} \label{dss5fd}
\mathbf{P}\bigg( \max_{1\leq i \leq n} |S_i|  \geq \lambda  \sqrt{n} \bigg) \ \leq \  \frac{ \varepsilon}{\lambda^{2}}.
\end{equation}
Since $$ \mathbf{E} [|\xi_i|^3 |\mathcal{F}_{i-1} ] \leq M  \mathbf{E} [\xi_i^2 |\mathcal{F}_{i-1} ],$$  we deduce that
$$ (\mathbf{E} [\xi_i^2 |\mathcal{F}_{i-1} ] )^{3/2} \leq  \mathbf{E} [|\xi_i|^3 |\mathcal{F}_{i-1} ]   \leq   M  \mathbf{E} [\xi_i^2 |\mathcal{F}_{i-1} ].$$
Thus $$ \mathbf{E} [\xi_i^2 |\mathcal{F}_{i-1} ]  \leq     M^2, \ \ \ \ \  \langle S\rangle_n \leq nM^2 \ \ \ \textrm{and} \  \ \ \sum_{i=1}^{n}\mathbf{E} [|\xi_i|^3 |\mathcal{F}_{i-1} ] \leq n M^3.$$
Applying (\ref{sdfds}) with $p=3, x=  2y=\lambda  \sqrt{n} $, we obtain
\begin{eqnarray}
\mathbf{P}\bigg( \max_{1\leq i \leq n} |S_i|  \geq \lambda  \sqrt{n} \bigg)  &\leq &    2\exp\Bigg\{- \frac{2 \, \lambda^2}{25e^3  M^2} \Bigg \} + 2\exp\Bigg\{- \frac{6}{5} \log\Bigg( 1 + \frac{3 \lambda^3 \sqrt{n} }{20 M^3}\Bigg)  \Bigg \} \ \nonumber  \\
&&\ \ \   + \ \mathbf{P}\left( \max_{ 1\leq i\leq n}\xi_{i} >  \frac12 \lambda  \sqrt{n}  \right) +  \mathbf{P}\left( \max_{ 1\leq i\leq n}(-\xi_{i}) > \frac12 \lambda  \sqrt{n}  \right)\, \nonumber\\
 &\leq &    2\exp\Bigg\{- \frac{4 \, \lambda^2}{50e^3  M^2} \Bigg \} + 2\Bigg(   \frac{3 \lambda^3 \sqrt{n} }{20 M^3}\Bigg)^{-6/5}  \ \nonumber  \\
&&\ \ \   + \ \frac{16}{\lambda^3 \sqrt{n}} \mathbf{E} \Big[|\xi_1|^3 \mathbf{1}_{\{|\xi_{1}| >  \frac12 \lambda  \sqrt{n}\}  } \Big]  \ \nonumber  \\
 &\leq &     \frac{ \varepsilon}{\lambda^{2}}
\end{eqnarray}
provided that $\lambda$ is sufficiently large. This proves (\ref{dss5fd}).\hfill\qed

\section*{Acknowledgements}
I gratefully think J\'{e}r\^{o}me Dedecker for his helpful discussion and corrections. I am indebted to
Haijuan Hu for her helpful corrections.

%I am grateful to the referees for their helpful remarks and suggestions.

\end{document}